\documentclass{article}  
\usepackage[pdftex,
bookmarksnumbered,
bookmarksopen,
colorlinks,
citecolor=blue,
linkcolor=blue,]{hyperref}
\usepackage{amsmath}
\usepackage{amssymb}
\usepackage{mathrsfs}
\usepackage{amsfonts}
\usepackage{amsthm}
\usepackage{epsfig}
\usepackage{graphicx}
\usepackage{subfigure}
\usepackage{epstopdf}
\usepackage{color}
\usepackage{algorithm}
\usepackage{booktabs}
\usepackage{listings}
\usepackage{boxedminipage}
\usepackage{algorithmic}
\usepackage{stmaryrd}
\usepackage{cite}
\usepackage{relsize}
\usepackage{lscape}
\usepackage[top=1.3in, bottom=1.5in, left=1.5in, right=1.5in]{geometry}
\usepackage{stmaryrd} 
\usepackage{relsize}%
\usepackage{url}
\usepackage{color,xcolor}
\usepackage{pgf,tikz,pgfplots}
\usepackage{multirow}

\newtheorem{theorem}{Theorem}[section]

\newtheorem{proposition}{Proposition}[section]

\newtheorem{assumption}{Assumption}[section]
\newtheorem{lemma}{Lemma}[section]
\newtheorem{corollary}{Corollary}[section]
\newtheorem{remark}{Remark}[section]

\newenvironment{mytabular}{\bgroup\tiny\tabular}{\endtabular\egroup}

\def\norm#1{\|#1\|}

\newcommand{\T}{\ensuremath{ \mathbb R^{n_1\times \cdots\times n_d}   }}

\newcommand{\bigxiaokuohao}[1]{\ensuremath{ \left(  #1 \right) }}      
\newcommand{\bigjueduizhi}[1]{\ensuremath{ \left|  #1 \right| }}   
\newcommand{\bigdakuohao}[1]{\ensuremath{ \left\{  #1 \right\} }}

\newcommand{\bigfnorm}[1]{\ensuremath{ \left\|   #1 \right\|_F }}    
\newcommand{\bignuclearnorm}[1]{\ensuremath{ \left\|   #1 \right\|_* }}    
\newcommand{\bignorm}[1]{\ensuremath{ \left\|   #1 \right\|  }}

\newcommand{\innerprod}[2]{\ensuremath{ \left\langle   #1 , #2\right\rangle }}      

\newcommand{\bigotimesu}{\ensuremath{     \bigotimes^d_{j=1}\nolimits \mathbf u_{j,i  } }}

\DeclareMathOperator*{\E}{\mathop{\mathsf{E}}}
\DeclareMathOperator*{\Einline}{\mathop{{}\mathsf{E}}}
                                                              
	\definecolor{darkgray}{rgb}{0.66, 0.66, 0.66}

\title{On Approximation Algorithm for Orthogonal Low-Rank Tensor Approximation
}

\author{Yuning Yang  \thanks{College of Mathematics and Information Science, Guangxi University, Nanning, 530004, China  (yyang@gxu.edu.cn).}                             
}
\begin{document} %\large
\maketitle

\begin{abstract}
The goal of this work is to fill a   gap   in [Yang, SIAM J. Matrix Anal. Appl, 41 (2020), 1797--1825]. In that work, an approximation procedure   was proposed for orthogonal low-rank tensor approximation; however, the approximation lower bound was only established when the number of orthonormal factors is one. To this end, by further  exploring the multilinearity and orthogonality of the problem, we introduce a modified approximation algorithm. Approximation lower bound is established, either in   deterministic or     expected sense, no matter how many orthonormal factors there are. In addition, a major feature of the new algorithm  is its flexibility to allow either deterministic or randomized procedures to solve a key step of each latent orthonormal factor  involved in the algorithm. This feature can reduce the computation of large   SVDs, making the algorithm more efficient.  Some numerical studies are provided to validate the usefulness of the proposed algorithm.

\end{abstract}

\section{Introduction}
Tensors (hypermatrices) play an important role in signal processing, data analysis, statistics, and machine learning nowadays, key to which are   tensor decompositions/approximations \cite{kolda2010tensor,comon2014tensors,cichocki2015tensor,Sidiropoulos2016ten}. 
Canonical polyadic (CP) decomposition factorizes the data tensor into some rank-1 components. Such a decomposition can be unique under mild assumptions \cite{kruskal1977three,qi2016semialgebraic,sidiropoulos2000uniqueness}; however, the degeneracy of its associated optimization requires one to impose additional constraints, such as the nonnegativity, angularity, and orthogonality \cite{kolda2001orthogonal,lim2013blind,lim2009nonnegative}. Orthogonal CP decomposition allows one or more of the latent factors to be orthonormal \cite{kolda2001orthogonal,guan2019numerical}, giving flexibility to model applications arising from image processing, joint SVD, independent component analysis,  DS-CDMA systems, and so on \cite{shashua2001linear,comon1994independent,sorensen2010parafac,pesquet2001joint,de2011short,sidiropoulos2000blind,sorensen2010parafac}. 

Most existing algorithms for solving orthogonal low-rank tensor decomposition/approximation are iterative algorithms  \cite{chen2009tensor,sorensen2012canonical,wang2015orthogonal,pan2018symmetric,guan2019numerical,yang2019epsilon,hu2019linear,li2019polar,li2018globally,savas2010quasi}; just to name a few. To find a good feasible initializer for iterative algorithms, \cite[Procedure 3.1]{yang2019epsilon} developed an approximation procedure, which finds the orthonormal factors by means of the celebrated   HOSVD \cite{de2000a}, and then computes the non-orthonormal ones by solving  tensor rank-1 approximation problems. Some numerical tests showed that with the help of this procedure,  the alternating least squares (ALS) can   recover the latent factors better.

However, a   gap was left: the approximation lower bound of    \cite[Procedure 3.1]{yang2019epsilon}  was established only when the number of orthonormal factors is one. 
Denote $t$ as the number of latent orthonormal factors and $d$ the order of the data tensor under consideration. The difficulty  to extend the lower bound analysis of  \cite[Procedure 3.1]{yang2019epsilon}  to arbitrary $1<t\leq d$ 
lies  in that the based HOSVD is essentially designed for the Tucker format, which would yield some cross-terms that cause troubles in the analysis of the current CP format.

To fill this   gap, in this work, by further exploring the multilinearity and orthogonality of the   problem, we devise a new approximation algorithm for orthogonal low-rank tensor CP approximation, which is a modification of \cite[Procedure 3.1]{yang2019epsilon}. In fact, we will make the new algorithm  somewhat   more general,  allowing either deterministic or randomized procedures to solve a key step of each latent orthonormal factor involved in the algorithm. Following a certain principle, three procedures for this step are considered, two deterministic and one randomized, leading to three versions of the approximation algorithm. 
The approximation lower bound is derived, either in the deterministic or in the expected sense, no matter how many latent orthonormal factors there are. The most expensive operation of the proposed algorithm is an   SVD of   the mode-$d$ unfolding matrix, making the algorithm more efficient than \cite[Procedure 3.1]{yang2019epsilon}.
Numerical studies will   show the efficiency of the introduced algorithm and its usefulness   in recovery and clustering, and is favorably compared with \cite[Procedure 3.1]{yang2019epsilon}.

Note that approximation algorithms have been widely proposed for tensor approximations.  For example, HOSVD, sequentially truncated HOSVD \cite{vannieuwenhoven2012new}, and hierarchical HOSVD \cite{grasedyck2010hierarchical} are also   approximation algorithms for the Tucker format. On the other hand, several randomized approximation algorithms have been proposed in recent years for Tucker format or t-SVD     \cite{minster2020randomized,che2020computation,zhang2018randomized,ahmadi2020randomized}. The solution quality of the aforementioned algorithms was measured via error bounds. In the context of tensor best rank-1 approximation problems, deterministic or randomized approximation algorithms were   developed \cite{he2010approximation,he2014probability,so2010deterministic,da2016finite,kofidis2002on}, where the solution quality was measured by approximation bounds; our results extend these to rank-$R$ approximations. 

The rest of this work is organized as follows. Preliminaries are given in Sect. \ref{sec:problem}. The approximation algorithm is introduced in Sect. \ref{sec:approx_alg}, while the approximation bound and its analysis are derived in Sect. \ref{sec:approx_bound}. Sect. \ref{sec:numer} provides numerical results and Sect. \ref{sec:conclusions} draws some conclusions.  

\section{Preliminaries} \label{sec:problem}
Throughout this work, vectors are written as boldface lowercase letters $(\mathbf x,\mathbf y,\ldots)$, matrices
correspond to italic capitals $(A,B,\ldots)$, and tensors are
written as calligraphic capitals $(\mathcal{A}, \mathcal{B},
\cdots)$. $\mathbb R^{n_1\times \cdots\times n_d}$ denotes the space of $n_1\times\cdots\times n_d$ real tensors. 
For two tensors $\mathcal A,\mathcal B$ of the same size, their inner product $\langle \mathcal A,\mathcal B\rangle$  is given by
the sum of entry-wise product. 
The Frobenius (or Hilbert–Schmidt)  norm of $\mathcal A$ is defined by $\|\mathcal A\|_F = \langle\mathcal A,\mathcal A\rangle^{1/2}$;    $\otimes$ denotes the outer product; in particular, for $\mathbf u_j\in\mathbb R^{n_j}$, $j=1,\ldots,d$, $\mathbf u_1\otimes\cdots\otimes\mathbf u_d$ denotes a rank-1 tensor in $\T$.  We write it as $\bigotimes^d_{j=1}\mathbf u_j$ for short. For  $R$ vectors of size $n_j$:  $\mathbf u_{j,1},\ldots,\mathbf u_{j,R} $, we usually collect them into a matrix $U_j=[\mathbf u_{j,1},\ldots,\mathbf u_{j,R}] \in\mathbb R^{n_j\times R}$, $j=1,\ldots,R$. 

The mode-$j$ unfolding of $\mathcal A\in\T$, denoted as $A_{(j)}$, is a matrix  in $\mathbb R^{n_j\times \prod^d_{k\neq j}n_k}$. 
The product between a tensor $\mathcal A\in \T$  and a vector $\mathbf u_j\in\mathbb R^{n_j}$ with respect to the $j$-th mode, the tensor-vector product $\mathcal A\times_j \mathbf u_j^\top$ is a tensor in $ \in \mathbb R^{n_1\times \cdots\times n_{j-1}\times n_{j+1}\times \cdots \times n_d}$.

Given a $d$-th order tensor $\mathcal A\in \T$ and a positive integer $R$, the   approximate CP decomposition   under consideration can be written as  \cite{guan2019numerical}:
 \begin{equation} \label{prob:orth_main}
 \begin{split}
 &\min  ~   \bigfnorm{ \mathcal A - \sum^R_{i=1}\nolimits\sigma_i\bigotimesu}^2  \\
 &~    {\rm s.t.}~ ~ \mathbf u_{j,i}^\top \mathbf u_{j,i} =1, j=1,\ldots,d-t,  1\leq i \leq R,\\
 & ~~~~~~~ U_j^\top U_j = I,  j=d-t+1,\ldots,d, ~ { \sigma}_i\in\mathbb R,
 \end{split}
 \end{equation}
 where $1\leq t\leq d$.
This model means that the last $t$ latent factor matrices of $\mathcal A$ are assumed to be orthonormal, while the first $d-t$ ones are not; however, due to the presence of the scalars $\sigma_i$, when $j=1,\ldots,d-t$, each $\mathbf u_{j,i}$ can be normalized. These lead to the two types of  constraints in \eqref{prob:orth_main}. Using orthogonality, \eqref{prob:orth_main} can be equivalently formulated as a maximization problem \cite{guan2019numerical}:
 \begin{equation}		
 	      \setlength\abovedisplayskip{3pt}
 \setlength\abovedisplayshortskip{3pt}
 \setlength\belowdisplayskip{3pt}
 \setlength\belowdisplayshortskip{3pt}
 \label{prob:ortho_main_max}
  \begin{split}
&\max~   G(U_1,\ldots,U_d):= \sum^R_{i=1}\nolimits  \innerprod{\mathcal A}{\bigotimesu}^2 \\
 &~    {\rm s.t.}~ ~ \mathbf u_{j,i}^\top \mathbf u_{j,i} =1, j=1,\ldots,d-t,  1\leq i \leq R,\\
& ~~~~~~~ U_j^\top U_j = I,  j=d-t+1,\ldots,d,
 \end{split}
 \end{equation}
 where the variables $\sigma_i$'s have been eliminated. 
The approximation algorithms and the corresponding analysis are mainly focused on  this maximization model.  
 
 Some necessary definitions and properties are introduced in the following. 
 
Let $V\in\mathbb R^{m\times n}$, $m\geq n$. The polar decomposition of $V$ is to decompose it into two matrices $U\in\mathbb R^{m\times n}$ and $H\in\mathbb R^{n\times n}$ such that $V=UH$, where $U$ is columnwisely orthonormal, and $H$ is a symmetric positive semidefinite matrix. Its relation with   SVD is given in the following. 
 \begin{theorem}[c.f. \cite{higham1986computing}] \label{th:polar_dec}
	Let $V\in\mathbb R^{m\times n}$, $m\geq n$. Then there exist 	 $U\in\mathbb R^{m\times n}$ and a unique symmetric positive semidefinite matrix $H\in\mathbb R^{n\times n}$ such that
	\[ \setlength\abovedisplayskip{3pt}
	\setlength\abovedisplayshortskip{3pt}
	\setlength\belowdisplayskip{3pt}
	\setlength\belowdisplayshortskip{3pt}
	V = UH,~U^\top U = I\in\mathbb R^{n\times n}.
	\]
	$(U,H)$ is the polar decomposition of $V$. If ${\rm rank}(V)=n$, then $H$ is symmetric positive definite and $U$ is uniquely determined.
	
	Furthermore, let $H=Q\Lambda Q^\top$, $Q,\Lambda\in\mathbb R^{n\times n}$ be the eigenvalue decomposition of $H$, namely, $Q^\top Q = QQ^\top = I$, $\Lambda = {\rm diag}(\lambda_1,\ldots,\lambda_n)$ be a diagonal matrix where $\lambda_1\geq \cdots\geq \lambda_n\geq 0$. Then   
	$
	U = PQ^\top$, and  $V = P\Lambda Q^\top $
	is   a reduced SVD of $V$.
\end{theorem}

Write $V=[\mathbf v_1,\ldots,\mathbf v_n] $ and $U=[\mathbf u_1,\ldots,\mathbf u_n]$ where $\mathbf v_i,\mathbf u_i\in\mathbb R^m$, $i=1,\ldots,n$. The following lemmas are useful in the approximation bound analysis.
\begin{lemma}
	\label{lem:polar_decomp}
Let $U\in\mathbb R^{m\times n}$ and $H\in\mathbb R^{n\times n}$ be the two matrices generated by the polar decomposition of $V\in\mathbb R^{m\times n}$ such that $V=UH$, where $U$ is columnwisely orthonormal and $H$ is symmetric positive semidefinite. 	Let $\mathbf v_i$ and $\mathbf u_i$ be defined as above. Then it holds that
\[	
\setlength\abovedisplayskip{3pt}
\setlength\abovedisplayshortskip{3pt}
\setlength\belowdisplayskip{3pt}
\setlength\belowdisplayshortskip{3pt}
\innerprod{\mathbf u_i}{\mathbf v_i} \geq 0,~i=1,\ldots,m.
\]
	\end{lemma}
\begin{proof}
	Since $V=UH$ and $U^\top U = I$, we have $U^\top V = H$, which means that $\innerprod{\mathbf u_i}{\mathbf v_i}$ is exactly the $i$-th diagonal entry of $H$. As $H$ is positive semidefinite, its diagonal entries are nonnegative. The required result follows.
	\end{proof} 

Denote $\bignuclearnorm{\cdot}$ as the nuclear norm of a matrix, i.e., the sum of its singular values.   From Theorem \ref{th:polar_dec} we easily see that:
\begin{lemma}\label{lem:polar_nuclear_norm}
	Let $U$ and $V$ be defined as in Theorem \ref{th:polar_dec}. Then $\innerprod{U}{V} = \bignuclearnorm{V}$. 
	\end{lemma}

\section{Approximation Algorithm}\label{sec:approx_alg}

 The approximation algorithm proposed in this section inherits certain ideas of  \cite[Procedure 3.1]{yang2019epsilon}, and so we briefly recall its strategy here. 
 \cite[Procedure 3.1]{yang2019epsilon}   obtained an approximation solution to \eqref{prob:ortho_main_max} as follows:  One first applies the truncated HOSVD \cite{de2000a} to get $U_d,\ldots,U_{d-t+1}$; specifically, for $j=d,\ldots,d-t+1$, one unfolds $\mathcal A$ to the matrix $A_{(j)} \in\mathbb R^{n_j\times \prod^d_{l\neq j}n_l}$, and then performs the truncated SVD to get its left leading $R$  singular vectors, denoted as $\mathbf u_{j,1},\ldots,\mathbf u_{j,R}$  , which forms the $j$-th factor $U_j=[\mathbf u_{j,1},\ldots,\mathbf u_{j,R}]$. Once $U_d,\ldots,U_{d-t+1}$ are obtained, one then approximately solves $R$ tensor best rank-1 approximation problems, given the data tensors  $\mathcal A\bigotimes^{d}_{j=d-t+1}\mathbf u_{j,i}:=\mathcal A\times_{d-t+1} \mathbf u_{d-t+1,i}^\top\times \cdots \times_d\mathbf u_{d,i}^\top\in\mathbb R^{n_1\times \cdots\times n_{d-t}}$, $i=1,\ldots,R$, namely,
\begin{equation}
\label{prob:rank1approx}		
\begin{split}
&\max~    \innerprod{ \mathcal A\bigotimes^{d}_{j=d-t+1}\nolimits\mathbf u_{j,i}  }{ \bigotimes^{d-t}_{j=1}\nolimits\mathbf u_{j,i}   }  \\
& ~~{\rm s.t.}~~  \mathbf u_{j,i}^\top \mathbf u_{j,i} =1, j=1,\ldots,d-t,
\end{split}
\end{equation}
for $ i=1,\ldots,R$. 

Such a strategy of obtaining an approximation solution to \eqref{prob:ortho_main_max} works well empirically, while the approximation bound was established only when $t=1$ \cite[Proposition 3.2]{yang2019epsilon}. We find that it is difficult to derive the approximation bound for general tensors when $t\geq 2$; this is because we unavoidably encounter the tensor $\mathcal A \times_{d-t+1}U_{d-t+1}^\top\cdots\times_d U_d^\top$ during the analysis. Such kind of tensors might be easier to deal with in Tucker model \cite{de2000a}; however, it consists of  cross-terms $\mathcal A\times_{d-t+1}\mathbf u_{d-t+1,i_{d-t+1}}^\top \cdots\times_d\mathbf u_{d,i_d}^\top$ with $i_{d-t+1}\neq \cdots\neq i_d$ that are not welcomed in problem \eqref{prob:ortho_main_max}.

 \begin{figure}[t]
	%\vspace{-0.3cm}
	\centering
	\begin{tikzpicture}[scale=1.2, transform shape]
	\pgfmathsetmacro{\cubex}{1.}
	\pgfmathsetmacro{\cubey}{1.}
	\pgfmathsetmacro{\cubez}{1.}
	\pgfmathsetmacro{\cubesmall}{1}
	\newcommand{\tensorx}{0.3}
	\newcommand{\svdx}{0.6}
	\newcommand{\svdxx}{\svdx+0.2}
	\newcommand{\svdy}{-0.5}
	\newcommand{\uthreex}{\svdx+1.1}
	\newcommand{\uthreexx}{\uthreex+0.6}
	\newcommand{\uthreey}{\svdy+0.2}
	\newcommand{\arrowxone}{\uthreex+1.25}
	\newcommand{\arrowyone}{\svdy+0.5}
	\newcommand{\tentimesuthreex}{\arrowxone+1.1}
	\newcommand{\tentimesuthreey}{\arrowyone+1.6}
	\newcommand{\arrowxtwo}{\tentimesuthreex+2.55}
	\newcommand{\arrowytwo}{\tentimesuthreey-0.2}
	\newcommand{\rectxone}{\arrowxtwo+1.1}
	\newcommand{\rectyone}{\arrowytwo+0.}
	\newcommand{\arrowxthree}{\rectxone+0.6}
	\newcommand{\arrowythree}{\rectyone}
	\newcommand{\rectxtwo}{\arrowxthree+0.7}
	\newcommand{\rectytwo}{\arrowythree}
	\newcommand{\arrowxfour}{\rectxtwo+0.2}
	\newcommand{\arrowyfour}{\rectytwo}
	\newcommand{\rectxthree}{\arrowxfour+1.65}
	\newcommand{\rectythree}{\svdy+0.2}
	\newcommand{\arrowxfive}{\rectxthree}
	\newcommand{\arrowyfive}{\rectythree-0.6}
	\newcommand{\rectxfour}{\arrowxfive}
	\newcommand{\rectyfour}{\arrowyfive-4.}
	\newcommand{\arrowxsix}{\arrowxfour+1}
	\newcommand{\arrowysix}{\rectyfour+0.3}
	\newcommand{\tentimesutwox}{\arrowxsix-3.5}
	\newcommand{\tentimesutwoy}{\arrowysix+1.6}
	\newcommand{\arrowxseven}{\tentimesutwox-1.1}
	\newcommand{\arrowyseven}{\tentimesutwoy-0.3}
	\newcommand{\rectxfive}{\arrowxseven-0.7}
	\newcommand{\rectyfive}{\arrowyseven}
	\newcommand{\arrowxeight}{\rectxfive-0.3}
	\newcommand{\arrowyeight}{\arrowysix}
	\newcommand{\rectxsix}{\arrowxeight-1.65}
	\newcommand{\rectysix}{\rectyfour}
	\newcommand{\arrowxnine}{\rectxsix-0.6}
	\newcommand{\arrowynine}{\rectysix}
	\newcommand{\rectxseven}{\arrowxnine-1.6}
	\newcommand{\rectyseven}{\arrowynine}
	\draw[black,fill=blue!20!white,line width=0.2ex] (\tensorx,0,0) -- ++(-\cubex,0,0) -- ++(0,-\cubey,0) -- ++(\cubex,0,0) -- cycle;
	\draw[black,fill=blue!30!white,line width=0.2ex] (\tensorx,0,0) -- ++(0,0,-\cubez) -- ++(0,-\cubey,0) -- ++(0,0,\cubez) -- cycle;
	\draw[black,fill=blue!30!white,line width=0.2ex] (\tensorx,0,0) -- ++(-\cubex,0,0) -- ++(0,0,-\cubez) -- ++(\cubex,0,0) -- cycle;
	
	\draw[->, line width=0.2ex] (\svdxx,\svdy+0.1) -- (\svdxx+0.9,\svdy+0.1);   
	\draw[color=black] (\svdxx+0.4,\svdy+0.3) node {SVD};
	
	%\node (rectangle) at (\uthreex,\uthreey) [draw,thick,minimum width=0.25cm,minimum height=1.cm, fill=blue!20] {};
	%\draw[color=black] (\uthreex,\uthreey-1.) node {$\footnotesize\mathbf u_{3,1}$};
	%
	%\draw[color=black] (\uthreex+0.5,\uthreey-0.15) node {$\cdots$};
	%
	%\node (rectangle) at (\uthreex+1.0,\uthreey) [draw,thick,minimum width=0.25cm,minimum height=1.cm, fill=blue!20] {};
	%\draw[color=black] (\uthreex+1.0,\uthreey-1.) node {$\footnotesize\mathbf u_{3,R}$};
	
	\node (rectangle) at (\uthreexx,\uthreey) [draw,thick,minimum width=1.cm,minimum height=1.cm, fill=red!35] {$U_3$};
	\draw[-,dashed, line width=0.1ex] (\uthreexx-0.25,\uthreey+0.5) -- (\uthreexx-0.25,\uthreey-0.5);   
	\draw[-,dashed, line width=0.1ex] (\uthreexx,\uthreey+0.5) -- (\uthreexx,\uthreey-0.5);   
	\draw[-,dashed, line width=0.1ex] (\uthreexx+0.25,\uthreey+0.5) -- (\uthreexx+0.25,\uthreey-0.5);   
	\draw[color=black] (\uthreexx-0.1,\uthreey-0.7) node {$\footnotesize\mathbf u_{3,1},\ldots,\mathbf u_{3,R}$};
	
	\draw[->, line width=0.2ex] (\arrowxone,\arrowyone) -- (\arrowxone+1.,\arrowyone+1.1);   
	\draw[->, line width=0.2ex] (\arrowxone,\arrowyone-0.2) -- (\arrowxone+1.,\arrowyone+0.1); 
	\draw[->, line width=0.2ex] (\arrowxone,\arrowyone-0.4) -- (\arrowxone+1.,\arrowyone-0.6);     
	\draw[->, line width=0.2ex] (\arrowxone,\arrowyone-0.6) -- (\arrowxone+1.,\arrowyone-1.5);     
	
	\draw[black,fill=blue!20!white,line width=0.2ex] (\tentimesuthreex,\tentimesuthreey,0) -- (\tentimesuthreex+\cubesmall,\tentimesuthreey,0) -- (\tentimesuthreex+\cubesmall,\tentimesuthreey-\cubesmall,0) -- (\tentimesuthreex,\tentimesuthreey-\cubesmall,0) -- cycle;
	\draw[black,fill=blue!30!white,line width=0.2ex] (\tentimesuthreex,\tentimesuthreey,0) -- 
	(\tentimesuthreex,\tentimesuthreey,-\cubesmall) -- (\tentimesuthreex+\cubesmall,\tentimesuthreey,-\cubesmall) --
	(\tentimesuthreex+\cubesmall,\tentimesuthreey,0)  --
	cycle;
	\draw[black,fill=blue!30!white,line width=0.2ex] (\tentimesuthreex+\cubesmall,\tentimesuthreey-\cubesmall,0) -- (\tentimesuthreex+\cubesmall,\tentimesuthreey-\cubesmall,-\cubesmall) -- 
	(\tentimesuthreex+\cubesmall,\tentimesuthreey,-\cubesmall) --
	(\tentimesuthreex+\cubesmall,\tentimesuthreey,0) --
	cycle;	
	
	\draw[color=black] (\tentimesuthreex+\cubesmall+1,\tentimesuthreey-0.2) node {$\footnotesize\times_3\mathbf u_{3,1}^\top$};  
	
	\draw[color=black] (\tentimesuthreex+1.,\tentimesuthreey-1.3) node {$\vdots$};
	
	\draw[color=black] (\tentimesuthreex+1.,\tentimesuthreey-2.2) node {$\vdots$};

	\draw[black,fill=blue!20!white,line width=0.2ex] (\tentimesuthreex,\tentimesuthreey-3,0) -- (\tentimesuthreex+\cubesmall,\tentimesuthreey-3,0) -- (\tentimesuthreex+\cubesmall,\tentimesuthreey-3-\cubesmall,0) -- (\tentimesuthreex,\tentimesuthreey-3-\cubesmall,0) -- cycle;
	\draw[black,fill=blue!30!white,line width=0.2ex] (\tentimesuthreex,\tentimesuthreey-3.0,0) -- 
	(\tentimesuthreex,\tentimesuthreey-3.0,-\cubesmall) -- (\tentimesuthreex+\cubesmall,\tentimesuthreey-3.0,-\cubesmall) --
	(\tentimesuthreex+\cubesmall,\tentimesuthreey-3.0,0)  --
	cycle;
	\draw[black,fill=blue!30!white,line width=0.2ex] (\tentimesuthreex+\cubesmall,\tentimesuthreey-3.0-\cubesmall,0) -- (\tentimesuthreex+\cubesmall,\tentimesuthreey-3.0-\cubesmall,-\cubesmall) -- 
	(\tentimesuthreex+\cubesmall,\tentimesuthreey-3.0,-\cubesmall) --
	(\tentimesuthreex+\cubesmall,\tentimesuthreey-3.0,0) --
	cycle;	
	
	\draw[color=black] (\tentimesuthreex+\cubesmall+1,\tentimesuthreey-3.2) node {$\footnotesize\times_3\mathbf u_{3,R}^\top$};  
	
	\draw[->, line width=0.2ex] (\arrowxtwo,\arrowytwo) -- (\arrowxtwo+0.5,\arrowytwo);     
	
	%\draw[color=black] (\arrowxtwo,\tentimesuthreey-1.4) node {$\vdots$};
	%
	%\draw[color=black] (\arrowxtwo,\tentimesuthreey-2.8) node {$\vdots$};
	
	\draw[->, line width=0.2ex] (\arrowxtwo,\arrowytwo-3.) -- (\arrowxtwo+0.5,\arrowytwo-3.);

	\node (rectangle) at (\rectxone,\rectyone) [draw,thick,minimum width=1.cm,minimum height=1.cm, fill=blue!20] {$M_{2,1}$};
	
	\draw[color=black] (\rectxone,\tentimesuthreey-1.3) node {$\vdots$};
	
	\draw[color=black] (\rectxone,\tentimesuthreey-2.2) node {$\vdots$};
	
	\node (rectangle) at (\rectxone,\rectyone-3.) [draw,thick,minimum width=1.cm,minimum height=1.cm, fill=blue!20] {$M_{2,R}$};
	
	\draw[->, line width=0.2ex] (\arrowxthree,\arrowythree) -- (\arrowxthree+0.5,\arrowythree);     
	
	%\draw[color=black] (\arrowxthree+0.2,\tentimesuthreey-1.4) node {$\vdots$};
	%
	%\draw[color=black] (\arrowxthree+0.2,\tentimesuthreey-2.8) node {$\vdots$};
	
	\draw[->, line width=0.2ex] (\arrowxthree,\arrowythree-3.) -- (\arrowxthree+0.5,\arrowythree-3.); 
	
	\node (rectangle) at (\rectxtwo,\rectyone) [draw,thick,minimum width=0.1cm,minimum height=1.cm, fill=blue!20] {};
	\draw[color=black] (\rectxtwo,\tentimesuthreey-0.9) node {$\mathbf v_{2,1}$};
	
	\draw[color=black] (\rectxtwo,\tentimesuthreey-1.3) node {$\vdots$};
	
	\draw[color=black] (\rectxtwo,\tentimesuthreey-2.2) node {$\vdots$};
	
	\node (rectangle) at (\rectxtwo,\rectyone-3) [draw,thick,minimum width=0.1cm,minimum height=1.cm, fill=blue!20] {};
	\draw[color=black] (\rectxtwo,\tentimesuthreey-3.9) node {$\mathbf v_{2,R}$};
	
	\draw[->, line width=0.2ex] (\arrowxfour,\arrowyone+1.1) -- (\arrowxfour+1.,\arrowyone+0);   
	\draw[->, line width=0.2ex] (\arrowxfour,\arrowyone+0.1) -- (\arrowxfour+1.,\arrowyone-0.2); 
	\draw[->, line width=0.2ex] (\arrowxfour,\arrowyone-0.6) -- (\arrowxfour+1.,\arrowyone-0.4);     
	\draw[->, line width=0.2ex] (\arrowxfour,\arrowyone-1.5) -- (\arrowxfour+1.,\arrowyone-0.6);     
	
	\node (rectangle) at (\rectxthree,\rectythree) [draw,thick,minimum width=1.cm,minimum height=1.cm, fill=blue!20] {$V_{2}$};
	
	\draw[->, line width=0.2ex] (\arrowxfive,\arrowyfive) -- (\arrowxfive,\rectyfour+0.6); 
	
	\draw[color=black] (\arrowxfive+0.3,\arrowyfive-1.6) node {PD};
	
	\node (rectangle) at (\rectxfour,\rectyfour) [draw,thick,minimum width=1.cm,minimum height=1.cm, fill=red!35] {$U_{2}$};
	\draw[-,dashed, line width=0.1ex] (\rectxfour-0.25,\rectyfour+0.5) -- (\rectxfour-0.25,\rectyfour-0.5);   
	\draw[-,dashed, line width=0.1ex] (\rectxfour,\rectyfour+0.5) -- (\rectxfour,\rectyfour-0.5);   
	\draw[-,dashed, line width=0.1ex] (\rectxfour+0.25,\rectyfour+0.5) -- (\rectxfour+0.25,\rectyfour-0.5);   
	\draw[color=black] (\rectxfour+0.15,\rectyfour-0.7) node {$\footnotesize\mathbf u_{2,1},\ldots,\mathbf u_{2,R}$};

	\draw[->, line width=0.2ex] (\arrowxsix,\arrowysix) -- (\arrowxsix-1.,\arrowysix+1.1);   
	\draw[->, line width=0.2ex] (\arrowxsix,\arrowysix-0.2) -- (\arrowxsix-1.,\arrowysix+0.1); 
	\draw[->, line width=0.2ex] (\arrowxsix,\arrowysix-0.4) -- (\arrowxsix-1.,\arrowysix-0.6);     
	\draw[->, line width=0.2ex] (\arrowxsix,\arrowysix-0.6) -- (\arrowxsix-1.,\arrowysix-1.5);     
	
	\draw[black,fill=blue!20!white,line width=0.2ex] (\tentimesutwox,\tentimesutwoy,0) -- ++(-\cubex,0,0) -- ++(0,-\cubey,0) -- ++(\cubex,0,0) -- cycle;
	\draw[black,fill=blue!30!white,line width=0.2ex] (\tentimesutwox,\tentimesutwoy,0) -- ++(0,0,-\cubez) -- ++(0,-\cubey,0) -- ++(0,0,\cubez) -- cycle;
	\draw[black,fill=blue!30!white,line width=0.2ex] (\tentimesutwox,\tentimesutwoy,0) -- ++(-\cubex,0,0) -- ++(0,0,-\cubez) -- ++(\cubex,0,0) -- cycle;
	\draw[color=black] (\tentimesutwox+\cubesmall+0.55,\tentimesutwoy-0.2) node {$\footnotesize\times_2\mathbf u_{2,1}^\top\times_3\mathbf u_{3,1}^\top$};  
	
	\draw[color=black] (\arrowxtwo+1.1,\tentimesutwoy-1.3) node {$\vdots$};
	
	\draw[color=black] (\arrowxtwo+1.1,\tentimesutwoy-2.2) node {$\vdots$};

	\draw[black,fill=blue!20!white,line width=0.2ex] (\tentimesutwox,\tentimesutwoy-3.,0) -- ++(-\cubex,0,0) -- ++(0,-\cubey,0) -- ++(\cubex,0,0) -- cycle;
	\draw[black,fill=blue!30!white,line width=0.2ex] (\tentimesutwox,\tentimesutwoy-3.,0) -- ++(0,0,-\cubez) -- ++(0,-\cubey,0) -- ++(0,0,\cubez) -- cycle;
	\draw[black,fill=blue!30!white,line width=0.2ex] (\tentimesutwox,\tentimesutwoy-3.,0) -- ++(-\cubex,0,0) -- ++(0,0,-\cubez) -- ++(\cubex,0,0) -- cycle;
	\draw[color=black] (\tentimesutwox+\cubesmall+0.6,\tentimesutwoy-3.15) node {$\footnotesize\times_2\mathbf u_{2,R}^\top\times_3\mathbf u_{3,R}^\top$};

	\draw[->, line width=0.2ex] (\arrowxseven,\arrowyseven) -- (\arrowxseven-0.5,\arrowyseven); 
	\draw[->, line width=0.2ex] (\arrowxseven,\arrowyseven-3.1) -- (\arrowxseven-0.5,\arrowyseven-3.1); 
	
	\node (rectangle) at (\rectxfive,\rectyfive) [draw,thick,minimum width=0.1cm,minimum height=1.cm, fill=blue!20] {};
	\draw[color=black] (\rectxfive,\rectyfive-0.7) node {$\mathbf v_{1,1}$};
	
	\draw[color=black] (\tentimesuthreex+1,\rectyfive-1.05) node {$\vdots$};
	
	\draw[color=black] (\tentimesuthreex+1,\rectyfive-1.9) node {$\vdots$};
	
	\node (rectangle) at (\rectxfive,\rectyfive-3) [draw,thick,minimum width=0.1cm,minimum height=1.cm, fill=blue!20] {};
	\draw[color=black] (\rectxfive,\rectyfive-3.7) node {$\mathbf v_{1,R}$};
	
	\draw[->, line width=0.2ex] (\arrowxeight,\arrowyeight+1.1) -- (\arrowxeight-1.,\arrowyeight);   
	\draw[->, line width=0.2ex] (\arrowxeight,\arrowyeight+0.1) -- (\arrowxeight-1.,\arrowyeight-0.2); 
	\draw[->, line width=0.2ex] (\arrowxeight,\arrowyeight-0.6) -- (\arrowxeight-1.,\arrowyeight-0.4);     
	\draw[->, line width=0.2ex] (\arrowxeight,\arrowyeight-1.5) -- (\arrowxeight-1.,\arrowyeight-0.6);     
	
	\node (rectangle) at (\rectxsix,\rectysix) [draw,thick,minimum width=1.cm,minimum height=1.cm, fill=blue!20] {$V_{1}$};
	
	\draw[->, line width=0.2ex] (\arrowxnine,\arrowynine) -- (\arrowxnine-1,\arrowynine);   
	\draw[color=black] (\arrowxnine-0.45,\arrowynine+0.2) node {PD};
	
	\node (rectangle) at (\rectxseven,\rectyseven) [draw,thick,minimum width=1.cm,minimum height=1.cm, fill=red!35] {$U_{1}$};
	\draw[-,dashed, line width=0.1ex] (\rectxseven-0.25,\rectyseven+0.5) -- (\rectxseven-0.25,\rectyseven-0.5);   
	\draw[-,dashed, line width=0.1ex] (\rectxseven,\rectyseven+0.5) -- (\rectxseven,\rectyseven-0.5);   
	\draw[-,dashed, line width=0.1ex] (\rectxseven+0.25,\rectyseven+0.5) -- (\rectxseven+0.25,\rectyseven-0.5);   
	\draw[color=black] (\rectxseven-0.1,\rectyseven-0.7) node {$\footnotesize\mathbf u_{1,1},\ldots,\mathbf u_{1,R}$};
	\end{tikzpicture}
	
	\caption{Workflow of the proposed Algorithm \ref{alg:main} for approximately solving  \eqref{prob:ortho_main_max} when $d=3$ and $t=3$. PD is short for polar decomposition.} 
	\label{fig:alg_d3t3}
\end{figure}
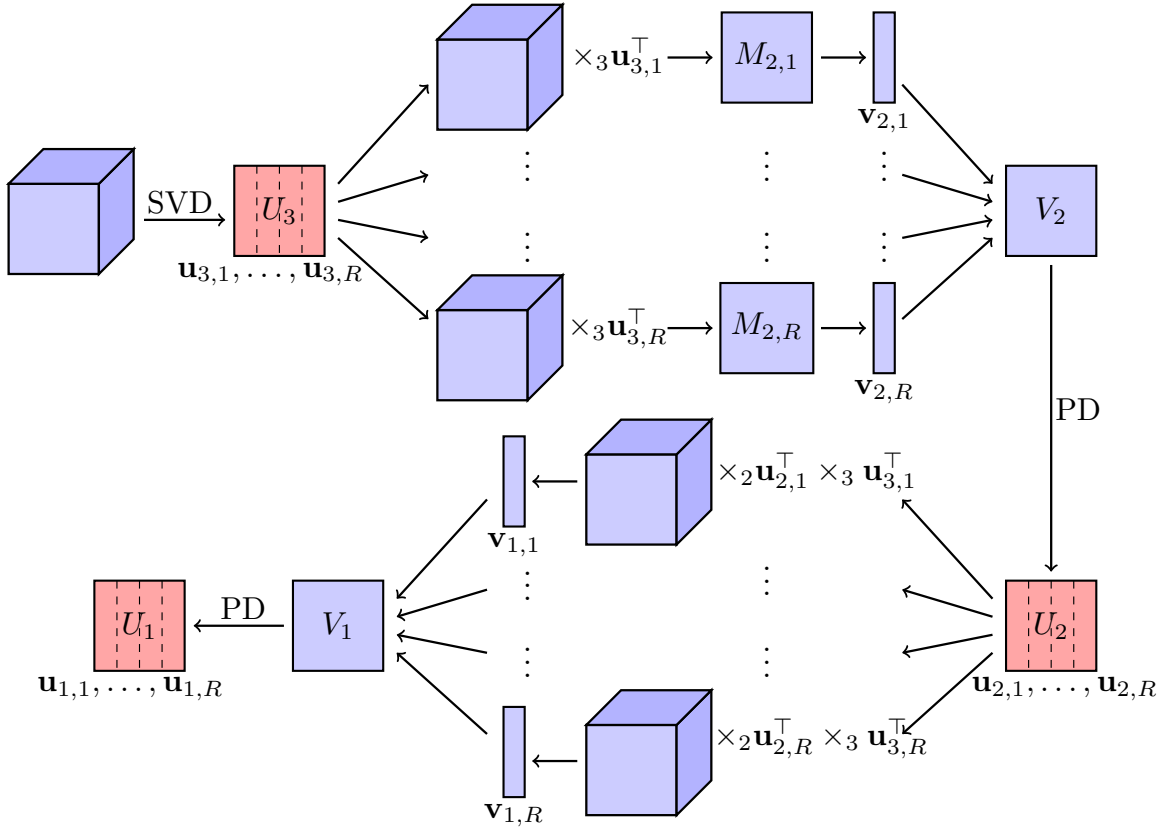

In view of this difficulty, we will modify \cite[Procedure 3.1]{yang2019epsilon} such that it is more suitable for       problem \eqref{prob:ortho_main_max}. 
The new algorithm is still designed under the guidance of problem \eqref{prob:rank1approx}, namely, we first compute the orthonormal factors $U_d,\ldots,U_{d-t+1}$, and then finding $U_1,\ldots,U_{d-t}$ via solving $R$  tensor rank-1 approximation problems. We differ from \cite{yang2019epsilon} in the computation of the orthonormal factors $U_{d-1},\ldots,U_{d-t+1}$.  By taking \eqref{prob:ortho_main_max} with $d= t=3$ as an example, we state our idea in what follows, with the workflow   illustrated in Fig. \ref{fig:alg_d3t3}. 

The computation of $U_3$ is the same as \cite{yang2019epsilon}, namely, we let $U_3:=[\mathbf u_{3,1},\ldots,\mathbf u_{3,R}]$ where $\mathbf u_{3,1},\ldots,\mathbf u_{3,R}$ are the left leading $R$ singular vectors of the unfolding matrix $A_{(d)} \in\mathbb R^{n_3\times n_1n_2}$. 

\begin{algorithm}[t]
	\algsetup{linenosize=\tiny}
	\footnotesize
	\caption{Approximation algorithm for   \eqref{prob:ortho_main_max}	\label{alg:main}}
	\begin{algorithmic}[1]
		\REQUIRE $\mathcal A\in\T$, $d\geq 3$, $R\geq 1$, $t\geq 1$ 
		\STATE Compute $\mathbf u_{d,1},\ldots,\mathbf u_{d,R}$ as the left leading $R$ unit singular vectors of the unfolding matrix $A_{(d)} \in\mathbb R^{n_d\times \prod^{d-1}_{k=1}n_k}$. Denote $U_d:= [\mathbf u_{d,1},\ldots,\mathbf u_{d,R}]$
		\FOR{$j = d-1:-1:d-t+1$} 
		\STATE \color{gray}\% splitting step\color{black}
		\FOR{$i = 1:R$}
		\STATE  Compute the $j$-th order tensor 
		\vspace*{-0.3\baselineskip} 
		$$\mathcal B_{j,i}:=\mathcal A \times_{j+1}\mathbf u_{j+1,i}^\top\times_{j+2} \cdots \times_d\mathbf u_{d,i}^\top \in\mathbb R^{n_1\times\cdots\times n_j} 
		\vspace*{-1\baselineskip} 
		$$
		\IF{$j>1$}  
		%\color{gray}\% {lines 4 and 5 can be done simultaneously for all $i=1,\ldots,R$}\color{black}
		\STATE Unfold  $\mathcal B_{j,i}$ to matrix $M_{j,i}$ as
		\vspace*{-0.5\baselineskip} 
		$$M_{j,i}:= \texttt{reshape}\left( \mathcal B_{j,i}, n_j, \prod^{j-1}_{k=1}\nolimits n_k \right) \in \mathbb R^{n_j \times \prod^{j-1}_{k=1} n_k}
		\vspace*{-0.7\baselineskip} 
		$$
		\STATE Extract a vector $\mathbf v_{j,i}$ from $M_{j,i}$ as
		$$~~~~~~~~~~ \mathbf v_{j,i}:= \texttt{get\_v\_from\_M}(M_{j,i})\in\mathbb R^{n_j} ~~\color{gray}\%  {\rm ~to~ be~ introduced~ later}\color{black}
		\vspace*{-1\baselineskip} 
		$$ 
		\ELSE%{~~~~~\color{gray}\% this case happens only when $t=d$}\color{black}
		\STATE  
		$~~~~~~~~~~~~~~~~~~~~~~\mathbf v_{j,i}:= \mathcal B_{j,i}\in\mathbb R^{n_j}$
		\ENDIF
		\ENDFOR
		\STATE  \color{gray}\% gathering step\color{black}
		\STATE				Denote $V_{j}:=[\mathbf v_{j,1},\ldots, \mathbf v_{j,R}]\in\mathbb R^{n_j\times R}$
		\STATE
		Compute the polar decomposition of $V_j$ to obtain orthonormal $U_j$
		$$\!\!\!\!\!\!\!\!\!\!\!\!\!\!\!\!\!\!\!\!\!\!\!\!\!\!\!\!\!\!\!\!\!\!\!\!\!\!\!\!\!\!\!\!\!\!\!\!\!\!\!\!\!\!\!\!\!\!\!\!\!\!\!\!U_j=[\mathbf u_{j,1},\ldots,\mathbf u_{j,R}]:= \texttt{polar\_decomp}(V_j)$$
		\ENDFOR{ {~~~~~~~~~~~~~~~~~~~~~~~~~~\color{gray}\% end of the computation of orthonormal factors }}\color{black}
		\IF{$t<d $}
		\FOR{$i = 1:R$}
		\STATE    Compute a rank-1 approximation solution to the tensor $\mathcal A\times_{d-t+1} \mathbf u_{d-t+1,i}^\top\times \cdots \times_d\mathbf u_{d,i}^\top = \mathcal B_{d-t+1,i}\times_{d-t+1}\mathbf u_{d-t+1,i}^\top\in \mathbb R^{n_1\times \cdots\times n_{d-t}}$, with each $\mathbf u_{j,i}$ normalized:
		\vspace*{-0.3\baselineskip} 
		\[
		(\mathbf u_{1,i},\ldots,\mathbf u_{d-t,i}) = \texttt{rank1approx}( \mathcal A\times_{d-t+1} \mathbf u_{d-t+1,i}^\top\times \cdots \times_d\mathbf u_{d,i}^\top   )
		\vspace*{-0.3\baselineskip} 
		\]
		\ENDFOR
		\STATE Collect 
		$
		U_j:=[\mathbf u_{j,1},\ldots,\mathbf u_{j,R}],~~j=1,\ldots,d-t
		$
		\ENDIF{~~~~~~~~~~~~~~~~~~~~~~~~~~\color{gray}\% end of the computation of non-orthonormal factors}\color{black} \color{black}
		\RETURN $U_1,\ldots,U_d$
	\end{algorithmic}
\end{algorithm}

Finding $U_2$ is a bit more complicated than \cite{yang2019epsilon}, which can be divided into the \emph{splitting} step and the \emph{gathering} step.
In the splitting step, with $\mathbf u_{3,1},\ldots,\mathbf u_{3,R}$ at hand, we first compute $R$ matrices $M_{2,1},\ldots,M_{2,R}$, with 
\[
M_{2,i}:= \mathcal A\times_3\mathbf u_{3,i}^\top \in\mathbb R^{n_2\times n_1},~i=1,\ldots,R.
\]
Then, using a procedure $\texttt{get\_v\_from\_M}$ that will be specified later, we compute $R$ vectors $\mathbf v_{2,i}$ of size $n_2$ from
$M_{2,i}$ for each $i$. The principle of $\texttt{get\_v\_from\_M}$ is to extract as much of the information from $M_{2,i}$ as possible, and to satisfy some theoretical bounds. We will consider three versions of such procedures, which will be detailed later. In the gathering step, we first denote $V_2:=[\mathbf v_{2,1},\ldots,\mathbf v_{2,R}  ] \in\mathbb R^{n_2\times R}$. As $V_2$ may not be orthogonal, it is not feasible. We thus propose to    apply  polar decomposition on $V_2$ to obtain the orthonormal matrix $U_2 = [\mathbf u_{2,1},\ldots,\mathbf u_{2,R}]$.

Computing $U_1$ is similar to that of $U_2$.  In the splitting step, with $U_3$ and $U_2$ at hand, we first compute $M_{1,i}:= \mathcal A\times_2\mathbf u_{2,i}^\top\times \mathbf u_{3,i}^\top\in\mathbb R^{n_1}$, $i=1,\ldots,R$. Since $M_{1,i}$'s are already vectors, we directly let $\mathbf v_{1,i}=M_{1,i}$. In the gathering step, we   let $V_1:=[\mathbf v_{1,1},\ldots,\mathbf v_{1,R}]$ and perform polar decomposition to get $U_1$. %The whole  workflow  is depicted in Fig. \ref{fig:alg_d3t3}.

For general $d\geq 3$, we can similarly apply the above splitting step and gathering step to obtain   $U_{d-1},U_{d-2},\ldots$ sequentially. In fact, the design of the splitting and gathering steps follows the principle that
\begin{equation}\label{eq:principle_inequality_tmp}
\cdots\geq \sum^R_{i=1}  \nolimits\innerprod{\mathbf u_{j,i}}{\mathbf v_{j,i}}^2 \geq \alpha \sum^R_{i=1}\nolimits\innerprod{\mathbf u_{j+1,i}}{\mathbf v_{j+1,i}}^2 \geq \cdots,
\end{equation}
which will be studied in details in the next section; here the factor $\alpha \in (0,1)$. If $t<d$, i.e., there exist non-orthonormal constraints in \eqref{prob:ortho_main_max}, then with     $U_{d-t+1},\ldots,U_d$ at hand, we then compute $U_1,\ldots,U_{d-t}$ via approximately solving $R$ tensor rank-1 approximation problems \eqref{prob:rank1approx}.
The whole algorithm is summarized in Algorithm \ref{alg:main}.

Some remarks on Algorithm \ref{alg:main} are given as follows.
\begin{remark}\label{rmk:alg}
	\begin{enumerate}
		\item When $t=1$, Algorithm \ref{alg:main} and \cite[Procedure 3.1]{yang2019epsilon} are exactly the same if they take the same rank-1 approximation procedure. 
		\item 	The computation of $\mathcal B_{j,i}$ can be more efficient: Since $\mathcal B_{j,i} = \mathcal A \times_{j+1}\mathbf u_{j+1,i}^\top\times_{j+2} \cdots \times_d\mathbf u_{d,i}^\top$ and $\mathcal B_{j-1,i} = \mathcal A \times_{j}\mathbf u_{j,i}^\top\times_{j+1} \cdots \times_d\mathbf u_{d,i}^\top$, we have the relation that 
		$$\mathcal B_{j-1,i} = \mathcal B_{j,i} \times_j\mathbf u_{j,i}^\top.$$ 
		This recursive definition of $\mathcal B_{j,i}$ reduces its computational complexity. 
		
		\item The splitting step can be computed in parallel. 
		
		\item If taking $R=1$ and $t=d$ in problem \eqref{prob:ortho_main_max}, we see that it is a best rank-1 approximation problem. Therefore, Algorithm \ref{alg:main} with $R=1$ and $t=d$ is capable of solving such a problem. 
	\end{enumerate}
\end{remark}

\paragraph{On the procedure $\texttt{get\_v\_from\_M}$} How to obtain $\mathbf v_{j,i}$ from $M_{j,i}$ is important both in theory and practice. %Recall that our principle is to keep as much of the information of $M_{j,i}$ as possible. To be more specific, 
To retain more information of $M_{j,i}$ in $\mathbf v_{j,i}$, we expect that 
\begin{equation} \bignorm{\mathbf v_{j,i}}^2 \geq \beta\bignorm{M_{j,i}}_F^2\label{eq:v_M_inequality_tmp}\end{equation}
for some factor $\beta\in(0,1)$. 
We thus provide three simple procedures to achieve this, two deterministic and one randomized. To simplify the notations we omit the footscripts $j$ and $i$.
The matrix $M$ in the procedures is of size $n\times m$, $\mathbf v\in\mathbb R^n$, and $\mathbf y\in\mathbb R^m$.

The first idea is straightforward: we let $\mathbf v$ be the leading left singular vector of $A$ times the leading singular value of $M$:

\begin{boxedminipage}{0.85\textwidth}\small
	\begin{equation}  \label{proc:get_v_from_M_A}
	\noindent {\rm Procedure}~ \mathbf v=  \texttt{get\_v\_from\_M}(M)  \tag{A}
	\end{equation}
1. Compute the leading right unit singular vector of $M$, denoted as $\mathbf y$;

2. Return $\mathbf v = M\mathbf y$. 
\end{boxedminipage}

The second one is to first pick the row of $M$ with the largest magnitude, and then normalize it to yield $\mathbf y$; finally let $\mathbf v=M\mathbf y$:

\begin{boxedminipage}{0.85\textwidth}\small
	\begin{equation}  \label{proc:get_v_from_M_B}
	\noindent {\rm Procedure}~ \mathbf v=  \texttt{get\_v\_from\_M}(M)  \tag{B}
	\end{equation}
	1. Denote $M^k$ as the $k$-th row of $M$. Let $M^{\bar k}$ be the row with the largest magnitude, i.e.,
	\[
	\bignorm{M^{\bar k}} = \max_{1\leq k\leq n}\bignorm{M^k};
	\]
	if there exist multiple $\bar k$ satisfying the above, choose the first one;
	
	2. 
	Denote $\mathbf y:= M^{\bar k}/\bignorm{M^{\bar k}}$;
	
	3. Return $\mathbf v = M\mathbf y$. 
\end{boxedminipage}

The third one is similar to the second one, but to be more efficient, we randomly and uniformly pick a row of $M$, normalize it to yield $\mathbf y$, and then let $\mathbf v = M\mathbf y$:

\begin{boxedminipage}{0.85\textwidth}\small
	\begin{equation}  \label{proc:get_v_from_M_C}
	\noindent {\rm Procedure}~ \mathbf v=  \texttt{get\_v\_from\_M}(M)  \tag{C}
	\end{equation}
	1. Randomly and uniformly choose $\bar k \in \{1,\ldots,n\}$; 
	
	2. Denote $\mathbf y:= M^{\bar k}/\bignorm{M^{\bar k}}$;
	
	3. Return $\mathbf v = M\mathbf y$. 
\end{boxedminipage}

We postpone the analysis of \eqref{eq:v_M_inequality_tmp} of each procedure in the next section. We also remark that $\mathbf y$ in each procedure is   important for the approximation bound analysis.

\paragraph{On the procedure $\texttt{rank1approx}$} To be more general,  the $\texttt{rank1approx}$ procedure  in Algorithm \ref{alg:main} can be any efficient approximation algorithm for solving tensor best rank-1 approximation problems, such as \cite[Algorithm 1  with DR 2]{he2010approximation},  \cite[Sect. 5]{zhang2012cubic},   \cite[Procedure 3.3]{yang2019epsilon}, or even Algorithm \ref{alg:main} itself (as that mentioned in Remark \ref{rmk:alg}). 
For the purpose of approximation bound analysis, we require $\texttt{rank1approx}$ to   have an approximation bound characterization, i.e., for any $m$-th order  data tensor $\mathcal C\in\mathbb R^{n_1\times\cdots\times n_m}$ (assuming that
$n_1\leq \cdots\leq n_m$ and
$m\geq 1$), the normalized solution $(\mathbf x_1,\ldots,\mathbf x_m)$ returned by $\texttt{rank1approx}$ admits an    approximation bound of the form:
\begin{equation}\label{eq:approx_bound_rank1approx} 
\begin{split}
&	 \innerprod{\mathcal C}{\bigotimes^m_{j=1}\nolimits\mathbf x_j}  
\geq     \frac{\bignorm{\mathcal C}_F} { \zeta(m)  },~\forall \mathcal C\in\mathbb R^{n_1\times\cdots\times n_m},~m\geq 1,
\end{split}
\end{equation}
where $1/\zeta(m) \leq 1$ represents the     factor. For \cite[Approximation 1 with DR 2]{he2010approximation}, it can be checked from \cite[Sect. 3.1]{he2010approximation} that when $m\geq 3$,
\[
\zeta(m) = \sqrt{ \prod^{m-2}_{j=1}\nolimits n_j  }\cdot\sqrt{n_1};
\]
for \cite[Procedure 3.3]{yang2019epsilon}, it follows from \cite[Proposition 3.1]{yang2019epsilon} that
\[
\zeta(m) =	\left\{\begin{array}{ll}
\sqrt{ \prod^{m-1}_{j=1}n_j\cdot \prod^{m/2-2}_{j=1}n_{2j+1} \cdot n_2^{-1}  }\sqrt{n_{m-1}n_m}, & m~{\rm even~and}~m\geq 4,\\
\sqrt{\prod^{m-1}_{j=2}n_j\cdot\prod^{(m+1)/2-2}_{j=1}n_{2j}}\sqrt{n_{m-1}n_m},& m~{\rm odd~and}~m\geq 3;
\end{array}\right.
\]
when $m=1,2$, namely, $\mathcal C$ is a vector or a matrix, we also have that $\zeta(2)=\sqrt{n_1}$ and $\zeta(1)=1$ for both the two algorithms.  

\paragraph{Computational complexity of Algorithm \ref{alg:main}}
 We analyze the computational complexity of Algorithm \ref{alg:main}. To simplify the presentation, we consider $n_1=\cdots=n_d=n$. Computing $U_d$ is    a truncated SVD, whose complexity is $O(n^2\cdot n^{d-1})$. In the computation of $U_{d-1}$, computing $\mathcal B_{d-1,i}$ and $M_{d-1,i}$ require $O(n^d)$ flops; since $M_{d-1,i}\in \mathbb R^{n\times n^{d-2}}$, computing $\mathbf v_{d-1,i}$ requires $O(n^2\cdot n^{d-2})$ flops if Procedure \ref{proc:get_v_from_M_A} is used, while it takes $O(n\cdot n^{d-2})$ flops for Procedures \ref{proc:get_v_from_M_B} and \ref{proc:get_v_from_M_C}. In any case, this is dominated by $O(n^d)$. Thus the complexity of the splitting step is $O(Rn^d )$. The main effort of the gathering step is the polar decomposition of $V_{d-1,i}\in \mathbb R^{R\times n}$, % which is essentially to compute the SVD of $V_{{d-1}}$ of size $n\times R$,  
  with complexity $O(R^2n)$. Hence computing $U_{d-1}$ requires 
$
 \label{eq:complexity_U_d-1}
 O(Rn^d) + O(R^2n)  
$
 flops. 
 
 The complexity of computing $U_{d-2}$ is similar: from Remark \ref{rmk:alg},  computing $\mathcal B_{d-2,i}$ and $M_{d-2,i}$ needs $O(n^{d-1})$ flops. The total complexity of $U_{d-2}$ is
$
 \label{eq:complexity_U_d-2}
  O(Rn^{d-1}) + O(R^2n) .
$

 For $U_{d-3},\ldots,U_{d-t+1}$ the analysis is analogous.  Denote $O(\texttt{rank1approx}(\cdot))$ as the complexity of the rank-1 approximation procedure, depending on which procedure is used. Then
 the updating of $U_1,\ldots,U_{d-t}$ has complexity 
 $$
  	      \setlength\abovedisplayskip{3pt}
 \setlength\abovedisplayshortskip{3pt}
 \setlength\belowdisplayskip{3pt}
 \setlength\belowdisplayshortskip{3pt}
 O(Rn^{d-t+1}) + O(R\cdot\texttt{rank1approx}(\mathcal B_{d-t+1,i}\times_{d-t+1}\mathbf u_{d-t+1,i}^\top)),$$
 where the first term comes from computing $\mathcal B_{d-t+1,i}$. 
  Summarizing the above discussions, we have
 \begin{proposition}[Computational complexity of Algorithm \ref{alg:main}]
 	\label{prop:complexity}
 	Assume that $n_1=n_2=\cdots=n_d=n$. The computational complexity of Algorithm \ref{alg:main} is
 	\[
 	 	      \setlength\abovedisplayskip{3pt}
 	\setlength\abovedisplayshortskip{3pt}
 	\setlength\belowdisplayskip{3pt}
 	\setlength\belowdisplayshortskip{3pt}
 O(n^{d+1}) + 	\sum^{t+1}_{j=2}   O(Rn^{d-j+2}) + O(tR^2n) +    O(R\cdot\texttt{rank1approx}(\mathcal B_{d-t+1,i}\times_{d-t+1}\mathbf u_{d-t+1,i}^\top));
 	\]
 	in particular, if $t=d$, then it is  
 	\[ 
 	 	      \setlength\abovedisplayskip{3pt}
 	\setlength\abovedisplayshortskip{3pt}
 	\setlength\belowdisplayskip{3pt}
 	\setlength\belowdisplayshortskip{3pt}
 	 	O(n^{d+1}) + 	\sum^{d+1}_{j=2} \nolimits  O(Rn^{d-j+2}) + O(dR^2 n) , 
 	\]
 	which is dominated by $O(n^{d+1})$. 
 	\end{proposition}
 
 Note that \cite[Procedure 3.1]{yang2019epsilon} requires $t$ SVDs of size $n\times n^d$, while Algorithm \ref{alg:main} performs only one SVD of this size plus additional operations of smaller size. This makes Algorithm \ref{alg:main} more efficient and will be confirmed by numerical observations.
  
 \section{Approximation Bound Analysis}\label{sec:approx_bound} This section is organized as follows: approximation bound results for general tensors are presented  in Sect. \ref{sec:approx_bound_general}, and then we shortly improve the results for nearly orthogonal tensors in Sect. \ref{sec:approx_bound_nearly_orth}.
 
 To begin with, we need some preparations.
 Recall that there is an intermediate variable $\mathbf y$ in Procedures \ref{proc:get_v_from_M_A}--\ref{proc:get_v_from_M_C}. This variable is important in the analysis. To distinguish $\mathbf y$ with respect to each $i$ and $j$, when obtaining $\mathbf v_{j,i}$, we denote the associated $\mathbf y$ as $\mathbf y_{j,i}$, wihch is of size $\prod^{j-1}_{k=1}n_k$. From the procedures, we see that $\mathbf v_{j,i} = M_{j,i}\mathbf y_{j,i}$. Since $M_{j,i} = \texttt{reshape}\left( \mathcal B_{j,i}, n_j, \prod^{j-1}_{k=1}\nolimits n_k \right) $ and  $\mathcal B_{j,i }=\mathcal A \times_{j+1}\mathbf u_{j+1,i}^\top\times_{j+2} \cdots \times_d\mathbf u_{d,i}^\top  $, we have the following expression of $\innerprod{\mathbf u_{j,i}}{\mathbf v_{j,i}}$:
 \begin{eqnarray}\label{eq:approx_bound_1}
 \innerprod{\mathbf u_{j,i}}{\mathbf v_{j,i}} &=& \innerprod{\mathbf u_{j,i}}{M_{j,i}\mathbf y_{j,i}}\nonumber\\
 & =& \innerprod{M^\top_{j,i}\mathbf  u_{j,i}}{\mathbf y_{j,i}} \\
 &=& \innerprod{\mathcal A\times_j\mathbf u_{j,i}^\top\times_{j+1}\cdots\times_d\mathbf u^\top_{d,i} }{\mathcal Y_{j,i}},\nonumber
 \end{eqnarray}
 where we denote
 \begin{equation}\label{eq:Yji}
 \mathcal Y_{j,i} := \texttt{reshape}(\mathbf y_{j,i},n_1,\ldots,n_{j-1}) \in\mathbb R^{n_1\times\cdots\times n_{j-1}}.
 \end{equation}
 Note that $\bigfnorm{\mathcal Y_{j,i}}=1$ because according to Procedures \ref{proc:get_v_from_M_A}--\ref{proc:get_v_from_M_C},  $\mathbf y_{j,i}$ is a unit vector. 
 
 The above expression is only well-defined when $j\geq 2$ (see Algorithm \ref{alg:main}). To make it consistent when $j=1$ (which happens when $t=d$), we    set $\mathcal Y_{1,i}=\mathbf y_{1,i}=1$; we also denote $M_{1,i}:=\mathcal B_{1,i} = \mathbf v_{1,i}$ accordingly.   
 Then it is clear that \eqref{eq:approx_bound_1} still makes sense when $j=1$.

 On the other hand, $V_j$ in the algorithm is only defined when $j=d-t+1,\ldots,d-1$. For convenience, we denote 
 \begin{equation}\label{eq:V_d}
 V_d=[\mathbf v_{d,1},\ldots,\mathbf v_{d,R}] \in\mathbb R^{n_d\times R} ~~{\rm with}~~ \mathbf v_{d,i}:= \lambda_i(A_{(d)}   )\mathbf u_{d,i}, ~i=1,\ldots, R, 
 \end{equation}
 where $\lambda_i(A_{(d)})$ denotes the $i$-th largest singular value of $A_{(d)}$.  
 
 Another notation to be defined is $\mathcal B_{d-t,i}$. Note that in the algorithm, $\mathcal B_{j,i}$ is only defined when $j\geq d-t+1$. We thus similarly define   $\mathcal B_{d-t,i}:= \mathcal A\times_{d-t+1}\mathbf u_{d-t+1,i}^\top\times_{d-t+2}\cdots\times_d\mathbf u_{d,i}^\top \in\mathbb R^{n_1\times\cdots\times n_{d-t}}$, $i=1,\ldots,R$. When $t=d$, $\mathcal B_{0,i}$'s are scalars.

 \subsection{Approximation bound for general tensors}\label{sec:approx_bound_general} 
   Since when $t=1$, Algorithm \ref{alg:main} coincides with \cite[Procedure 3.1]{yang2019epsilon} if they take the same best rank-1 approximation procedure,   similar to \cite[Proposition 3.2]{yang2019epsilon}, we   have that when $t=1$: 
   \begin{small}
 \begin{equation}
 \label{eq:tis1_case}
 G(U_1,\ldots,U_d) = \sum^R_{i=1}  \innerprod{\mathcal A\times_d\mathbf u_{d,i}^\top}{\bigotimes^{d-1}_{j=1} \mathbf u_{j,i}}^2 \geq \sum^R_{i=1}\frac{\bigfnorm{\mathcal A\times_d\mathbf u^\top_{d,i}}^2}{\zeta(d-1)^2} = \sum^R_{i=1}\frac{\lambda_i(A_{(d)} )^2}{\zeta(d-1)^2} ,
 \end{equation}
 \end{small}
 where the inequality comes from    \eqref{eq:approx_bound_rank1approx}, and the last equality is due to that the unfolding of $\mathcal A\times_d\mathbf u_{d,i}^\top$ to a vector is exactly $\mathbf v_{d,i}$ defined in \eqref{eq:V_d}. 
 
 Thus the remaining part is mainly focused on $t>2$ cases. To make the analysis more readable, 
 we first present an overview and informal analysis of the approximation bound, under the setting that \eqref{eq:principle_inequality_tmp} holds with a  factor $\alpha_j$, and then we prove the validity of \eqref{eq:principle_inequality_tmp} and detail the associated factor.   The formal approximation bound results are stated in the last of this subsection.

 \begin{lemma}[Informal approximation bound of Algorithm \ref{alg:main}]
 	\label{lem:approx_bound_informal_alg}
 	Let $2\leq t\leq d$ and let $U_1,\ldots,U_d$ be generated by Algorithm \ref{alg:main}. If for each $d-t+1\leq j\leq d-1$,  there holds
 	\begin{equation}
 	\label{eq:principle_inequality_tmp_alpha_j}
 	   \sum^R_{i=1}  \nolimits\innerprod{\mathbf u_{j,i}}{\mathbf v_{j,i}}^2 \geq \alpha_j \sum^R_{i=1}\nolimits \innerprod{\mathbf u_{j+1,i}}{\mathbf v_{j+1,i}}^2,
 	\end{equation}
 	\end{lemma}
 where $\alpha_j\in(0,1]$, then we have the following approximation bound:
 \begin{eqnarray*}
 \label{eq:approx_bound_informal}
 G(U_1,\ldots,U_d) &\geq & \frac{1}{\zeta(d-t)^2}\prod^{d-1}_{j=d-t+1}\nolimits\alpha_j \cdot \sum^R_{i=1} \nolimits \lambda_i(A_{(d)})^2\\
 & \geq& \frac{1}{ \zeta(d-t)^2}  \prod^{d-1}_{j=d-t+1}\nolimits\alpha_j \cdot G_{\max},
 \end{eqnarray*}
where $1/\zeta(d-t)$ is the rank-1 approximation ratio of  $n_1\times\cdots\times n_{d-t}$ tensors, as noted in \eqref{eq:approx_bound_rank1approx}.   
\begin{proof}
To make the analysis below consistent with the case $t=d$, which does not require perform rank-1 approximation,       we   denote $\zeta(0)=1$ and $\bigotimes^0_{j=1}\mathbf u_{j,i} = 1$.   It then holds that
	\begin{eqnarray*}
G(U_1,\ldots,U_d) &=& \sum^R_{i=1}\innerprod{\mathcal A}{\bigotimesu}^2 \\
&=& \sum^R_{i=1} \innerprod{\mathcal B_{d-t,i}}{ \bigotimes^{d-t}_{j=1}\nolimits \mathbf u_{j,i} }^2\\
&\geq & \frac{1}{\zeta(d-t)^2} \sum^R_{i=1} \bigfnorm{\mathcal B_{d-t,i}}^2 \\
&\geq& \frac{1}{\zeta(d-t)^2} \sum^R_{i=1} \innerprod{\mathcal B_{d-t,i}}{\mathcal Y_{d-t+1,i}}^2,
	\end{eqnarray*}
where the first inequality follows from the setting \eqref{eq:approx_bound_rank1approx}, and the second one is due to that $\bigfnorm{\mathcal B_{d-t,i}} = \max_{\bigfnorm{\mathcal Y}=1}\innerprod{\mathcal B_{d-t,i}}{\mathcal Y}$, and that $\bigfnorm{\mathcal Y_{d-t+1,i}}=1$  and $\mathcal Y_{d-t+1,i}\in\mathbb R^{n_1\times\cdots\times n_{d-t}}$  has the same size as $\mathcal B_{d-t,i}$; recall that $\mathcal Y_{j,i}$ was defined in \eqref{eq:Yji}.   From the definition of $\mathcal B_{d-t,i}$ and \eqref{eq:approx_bound_1}, we have
\begin{eqnarray*}
\innerprod{\mathcal B_{d-t,i}}{\mathcal Y_{d-t+1,i}} &=& \innerprod{\mathcal A\times_{d-t+1}\mathbf u_{d-t+1,i}^\top\times_{d-t+2}\cdots\times_d\mathbf u_{d,i}^\top}{\mathcal Y_{d-t+1,i}}\\
&=& \innerprod{\mathbf u_{d-t+1,i}}{\mathbf v_{d-t+1,i}}.
\end{eqnarray*}
It thus  follows from \eqref{eq:principle_inequality_tmp_alpha_j} that
\begin{eqnarray*}
G(U_1,\ldots,U_d) &\geq & \frac{1}{\zeta(d-t)^2}\sum^R_{i=1} \innerprod{\mathbf u_{d-t+1,i}}{\mathbf v_{d-t+1,i}}^2 \\
&\geq&\cdots\\
&\geq& \frac{1}{\zeta(d-1)^2}\prod^{d-t}_{j=d-t+1}\nolimits\alpha_j \cdot  \sum^R_{i=1}\nolimits \innerprod{\mathbf u_{d,i}}{\mathbf v_{d,i}}^2.
\end{eqnarray*}
From \eqref{eq:V_d} and that $\bignorm{\mathbf u_{j,i}}=1$, we see that
$
\innerprod{\mathbf u_{j,i}}{\mathbf v_{j,i}}^2 = \sum^R_{i=1}\nolimits\lambda_i(A_{(d)})^2,
$
which is clearly an upper bound of $G_{\max}$. Thus the required approximation bound follows.
	\end{proof}

The above lemma gives an overview of the analysis of the approximation bound. It also shows that to make the analysis go through, the chain inequality \eqref{eq:principle_inequality_tmp_alpha_j}  plays a central role. Thus the remaining task of this section is to establish \eqref{eq:principle_inequality_tmp_alpha_j} with an explicit factor $\alpha_j$.

 \subsubsection{Chain inequality \eqref{eq:principle_inequality_tmp_alpha_j} analysis}
To establish the connection between $\innerprod{\mathbf u_{j,i}}{\mathbf v_{j,i}}$ and $\innerprod{\mathbf u_{j+1,i}}{\mathbf v_{j+1,i}}$, we need to evaluate the performance of \texttt{get\_v\_form\_M}, namely, we need to estimate \eqref{eq:v_M_inequality_tmp}. As there are three choices of \texttt{get\_v\_form\_M}, we first simply assume that \eqref{eq:v_M_inequality_tmp} holds with a factor $\beta_j$, and later we analyze it in details. 
\begin{lemma}[Informal chain inequality]
\label{lem:approx_bound_informal_chain_inequality}
	Let $2\leq t\leq d$ and let $U_1,\ldots,U_d$ be generated by Algorithm \ref{alg:main}. If for each $d-t+1\leq j\leq d-1$, there holds
	\begin{equation}
	\label{eq:v_M_inequality_tmp_j}
	\bignorm{\mathbf v_{j,i}}^2 \geq \beta_j \bigfnorm{M_{j,i}}^2, ~i=1,\ldots,R,
	\end{equation}
	where $\beta_j\in(0,1]$ (if $j=1$ then \eqref{eq:v_M_inequality_tmp_j} holds with $\beta_j=1$), then we have for $j=d-t+1,\ldots,d-1$, and $i=1,\ldots,R,$
	\[
	\sum^R_{i=1} \innerprod{\mathbf u_{j,i}}{\mathbf v_{j,i}}^2 \geq \frac{\beta_j}{R}\sum^R_{i=1}\innerprod{\mathbf u_{j+1,i}}{\mathbf v_{j+1,i}}^2.
	\]
	\end{lemma}
\begin{remark}
	When $j=1$ (this happens when $t=d$), according to Algorithm \ref{alg:main}, $\mathbf v_{1,i}=\mathcal B_{1,i}=M_{1,i}$ for each $i$, and so \eqref{eq:v_M_inequality_tmp_j} holds with $\beta_1=1$.
\end{remark}
\begin{proof}
	Since the orthonormal $U_j$'s, $j=d-t+1,\ldots,d-1$, are given by the polar decomposition of $V_j$, from Lemma \ref{lem:polar_decomp}, we have $\innerprod{\mathbf u_{j,i}}{\mathbf v_{j,i}}\geq 0$. Using the relation that $a_1^2 + \cdots + a_R^2 \geq \frac{1}{R}\bigxiaokuohao{a_1+\cdots + a_R}^2$ for nonnegative $a_i$'s, we have 
	\[
	\sum^R_{i=1}\nolimits \innerprod{\mathbf u_{j,i}}{\mathbf v_{j,i}}^2 \geq \frac{1}{R}\bigxiaokuohao{\sum^R_{i=1}\nolimits\innerprod{\mathbf u_{j,i}}{\mathbf v_{j,i}} }^2 = \frac{1}{R}  \innerprod{U_j}{V_j}^2
	\]
for each $j$. 	Lemma \ref{lem:polar_nuclear_norm} shows that $\innerprod{U_j}{V_j} = \bignuclearnorm{V_j}$; the norm relationship gives that $\bignuclearnorm{V_j}\geq \bigfnorm{V_j}$. Thus
	\[
	\sum^R_{i=1}\nolimits \innerprod{\mathbf u_{j,i}}{\mathbf v_{j,i}}^2 \geq \frac{1}{R}\bigfnorm{V_j}^2.
	\]
	 We distinguish the cases that $j\leq d-2$ and $j=d-1$. For the former case, 
	\begin{eqnarray}
		\sum^R_{i=1}\nolimits \innerprod{\mathbf u_{j,i}}{\mathbf v_{j,i}}^2 &\geq& \frac{1}{R}\bigfnorm{V_j}^2 =\frac{1}{R}\sum^R_{i=1}\nolimits\bignorm{\mathbf v_{j,i}}^2\nonumber\\
	\footnotesize{{\rm (by~\eqref{eq:v_M_inequality_tmp_j})}}~~~~	 &\geq& \frac{\beta_j}{R} \sum^R_{i=1}\nolimits\bigfnorm{M_{j,i}}^2 \nonumber\\
		 &=& \frac{\beta_j}{R}\sum^R_{i=1}\nolimits\bigfnorm{\mathcal A\times_{j+1}\mathbf u^\top_{j+1,i}\times_{j+2}\cdots\times_d\mathbf u^\top_{d,i}  }^2 \label{eq:approx_bound_2}\\
		 &\geq& \frac{\beta_j}{R}\sum^R_{i=1}\nolimits\innerprod{ \mathcal A\times_{j+1}\mathbf u^\top_{j+1,i}\times_{j+2}\cdots\times_d\mathbf u^\top_{d,i}  }{\mathcal Y_{j+1,i}}^2\nonumber\\
\footnotesize{{\rm (by~\eqref{eq:approx_bound_1})}}	~~~~	 &=&\frac{\beta_j}{R} \sum^R_{i=1}\nolimits\innerprod{\mathbf u_{j+1,i}}{\mathbf v_{j+1,i}}^2,\nonumber
	\end{eqnarray}
	where the equality in the third line follows from that $M_{j,i}$ is the unfolding of $\mathcal B_{j,i} = \mathcal A \times_{j+1}\mathbf u^\top_{j+1,i}\times_{j+2}\cdots\times_d\mathbf u^\top_d$, and the third inequality uses the fact that $\bigfnorm{\mathcal Y_{j+1,i}}  =1$ where $\mathcal Y_{j,i}$ was defined in \eqref{eq:Yji}. 
	
	When $j=d-1$, since $\mathcal Y_{d,i}$ is not defined, there is a slight difference after the third line of \eqref{eq:approx_bound_2}. In this case, we have
	\begin{eqnarray*}
	\sum^R_{i=1}\nolimits\innerprod{\mathbf u_{d-1,i}}{\mathbf v_{d-1,i}}^2 &\geq& \frac{\beta_{d-1}}{R}\sum^R_{i=1}\nolimits\bigfnorm{ \mathcal A\times_d\mathbf u^\top_{d,i}   }^2\\
	& =& \frac{\beta_{d-1}}{R}\sum^R_{i=1}\nolimits\lambda_i(A_{(d)})^2 = \frac{\beta_{d-1}  }{R}\sum^R_{i=1}\nolimits\innerprod{\mathbf u_{d,i}}{\mathbf v_{d,i}}^2,
	\end{eqnarray*}
where the first inequality comes from a similar analysis as \eqref{eq:approx_bound_2}, the first equality follows from the definition of $U_d$,  and the last one is due to the definition of $V_d$ in \eqref{eq:V_d}. 
	The required inequality thus follows.
	\end{proof}

The remaining task is to specify \eqref{eq:v_M_inequality_tmp_j}.  
%Before this we need the following properties.  
Denote $\mathop{{}\mathsf{E}}$ as the expectation operator. 
\begin{lemma}
	\label{lem:approx_bound_get_v_from_M}
For $2\leq j\leq d-1$,	if $\mathbf v_{j,i}$ is generated from $M_{j,i}$  by Procedures \ref{proc:get_v_from_M_A} or \ref{proc:get_v_from_M_B}, then it holds that
	\begin{equation*}
\bignorm{\mathbf v_{j,i}}^2 \geq \frac{1}{n_j} \bigfnorm{M_{j,i}}^2;
	\end{equation*}
	if $\mathbf v_{j,i}$ is generated by Procedure  \ref{proc:get_v_from_M_C}, then it holds that
	\begin{equation*}
\E\bignorm{\mathbf v_{j,i}}^2 \geq \frac{1}{n_j}  \bigfnorm{M_{j,i}}^2.
	\end{equation*}
	\end{lemma}
\begin{proof}
	If $\mathbf v_{j,i}$ is generated by Procedure \ref{proc:get_v_from_M_A}, we have 
	\[
	\bignorm{\mathbf v_{j,i}}^2 =\lambda_{\max}^2(M_{j,i} ) \geq \frac{1}{n_j} \bigfnorm{M_{j,i}}^2, 
	\]
	where $\lambda_{\max}(\cdot)$ denotes the largest singular value of a matrix.
	
	If $\mathbf v_{j,i}$ is generated by Procedure \ref{proc:get_v_from_M_B}, we have
	\[
	\bignorm{\mathbf v_{j,i}}^2 = \bignorm{M_{j,i}\mathbf y_{j,i}}^2 = \sum^{n_j}_{k=1}  \innerprod{M_{j,i}^k}{\frac{M_{j,i}^{\bar k}  }{\bignorm{M_{j,i}^{\bar k}}}   }^2 \geq \bignorm{M_{j,i}^{\bar k}}^2 \geq \frac{1}{n_j}\bigfnorm{M_{j,i}}^2,
	\]
	where     we recall that $M_{j,i}^k$ denotes the $k$-th row of $M_{j,i}$, and $M_{j,i}^{\bar k}$ is the row with the largest magnitude. 
	
	If $\mathbf v_{j,i}$ is generated by Procedure \ref{proc:get_v_from_M_C}, since $\bar k$ is uniformly and randomly chosen from $\{1,\ldots,n_j\}$, this means that with equal probability, the random variable $\bignorm{\mathbf v_{j,i}}^2$ takes the value $\sum^{n_j}_{k=1}\innerprod{ M_{j,i}^k  }{ M_{j,i}^s/\bignorm{M_{j,i}^s}    }^2$, $s=1,\ldots,n_j$, i.e.,
	\[
	\mathop{\mathsf{Prob}}\bigdakuohao{\bignorm{\mathbf v_{j,i}}^2 =  \sum^{n_j}_{k=1} \innerprod{M_{j,i}^k}{\frac{M_{j,i}^s  }{\bignorm{M_{j,i}^s}}   }^2 } = \frac{1}{n_j},~s=1,\ldots,n_j.
	\] 
	Therefore,
	\begin{equation*}
\E \bignorm{\mathbf v_{j,i}}^2  =  \frac{1}{n_j}\sum^{n_j}_{s=1} \sum^{n_j}_{k=1} \innerprod{M_{j,i}^k}{\frac{M_{j,i}^s  }{\bignorm{M_{j,i}^s}}   }^2  
 \geq   \frac{1}{n_j} \sum^{n_j}_{s=1} \bignorm{M_{j,i}^s}^2 =\frac{1}{n_j}\bigfnorm{M_{j,i}  }^2. 
	\end{equation*}
The proof has been completed. 
	\end{proof}
\subsubsection{Puting the pieces together}
Combining the previous analysis, we can state our formal results.
\begin{lemma}[Formal chain inequality]
	\label{th:chain_inequality_formal}
		Let $2\leq t\leq d$ and let $U_1,\ldots,U_d$ be generated by Algorithm \ref{alg:main}. Let $\beta_j$ be such that
		\[
		\beta_j = {n_j^{-1}} ~{\rm if ~} 2\leq j\leq d-1,~~{\rm and}~~ \beta_j = 1 {\rm~if}~j=1.
		\]
		  For $j=d-t+1,\ldots,d-1,~i=1,\ldots,R$, if $\mathbf v_{j,i}$'s are generated deterministically by Procedures \ref{proc:get_v_from_M_A} or \ref{proc:get_v_from_M_B}, then
	\[
	\sum^R_{i=1}\nolimits \innerprod{\mathbf u_{j,i}}{\mathbf v_{j,i}}^2 \geq \frac{\beta_j}{R}\sum^R_{i=1}\nolimits\innerprod{\mathbf u_{j+1,i}}{\mathbf v_{j+1,i}}^2;
	\]
	if $\mathbf v_{j,i}$'s are generated randomly by Procedure  \ref{proc:get_v_from_M_C}, then
	\begin{equation*}
	\E \sum^R_{i=1}\nolimits \innerprod{\mathbf u_{j,i}}{\mathbf v_{j,i}}^2 \geq \frac{\beta_j}{R}\E \sum^R_{i=1}\nolimits\innerprod{\mathbf u_{j+1,i}}{\mathbf v_{j+1,i}}^2. 
	\end{equation*}
\end{lemma}
\begin{proof}
  Lemma \ref{lem:approx_bound_get_v_from_M} tells us that the factor $\beta_j$ in the assumption of  Lemma \ref{lem:approx_bound_informal_chain_inequality} is $ {\eta(j)}$. Thus	the results follow by combining Lemmas \ref{lem:approx_bound_informal_chain_inequality} and \ref{lem:approx_bound_get_v_from_M}. 
	\end{proof}
\begin{remark}
	The expectation here is taken sequentially with respect to $j=d-1,\ldots,d-t+1$, i.e., it is in fact 
	 \[ \E_{V_{d-1}}\nolimits\cdots \E_{V_j}\nolimits \sum^R_{i=1}\nolimits \innerprod{\mathbf u_{j,i}}{\mathbf v_{j,i}}^2  \geq \frac{\beta_j}{R}\E_{V_{d-1}}\nolimits\cdots \E_{V_{j+1}}\nolimits \sum^R_{i=1}\nolimits\innerprod{\mathbf u_{j+1,i}}{\mathbf v_{j+1,i}}^2.\]
\end{remark}

Based on \eqref{eq:tis1_case}, Lemmas \ref{lem:approx_bound_informal_alg} and \ref{th:chain_inequality_formal}, the main results are stated as follow.
\begin{theorem}[Formal approximation bound]
	\label{th:approx_bound}
	Let $1\leq t\leq d$ and let $U_1,\ldots,U_d$ be generated by Algorithm \ref{alg:main}. 	Let $\beta_j$ be defined as in Lemma \ref{th:chain_inequality_formal}.  For $j=d-t+1,\ldots,d-1,~i=1,\ldots,R$, if $\mathbf v_{j,i}$'s are generated deterministically by Procedures \ref{proc:get_v_from_M_A} or \ref{proc:get_v_from_M_B}, then 
 \begin{eqnarray*}
	G(U_1,\ldots,U_d) &\geq & \frac{1}{\zeta(d-t)^2}  \prod^{d-1}_{j=d-t+1}\nolimits\frac{\beta_j}{R} \cdot \sum^R_{i=1}\nolimits \lambda_i(A_{(d)})^2\\
	& \geq&  \frac{1}{\zeta(d-t)^2 } \prod^{d-1}_{j=d-t+1}\nolimits\frac{\beta_j}{R}  \cdot G_{\max};
\end{eqnarray*}
if 	  $\mathbf v_{j,i}$'s are generated randomly by Procedure  \ref{proc:get_v_from_M_C}, then
 \begin{eqnarray*}
	\Einline G(U_1,\ldots,U_d) &\geq &\frac{1}{ \zeta(d-t)^2 }\prod^{d-1}_{j=d-t+1}\nolimits\frac{\beta_j}{R}  \cdot\sum^R_{i=1} \nolimits \lambda_i(A_{(d)})^2\\
	& \geq& \frac{1}{ \zeta(d-t)^2 } \prod^{d-1}_{j=d-t+1}\nolimits\frac{\beta_j}{R}  \cdot G_{\max}. 
\end{eqnarray*}
The ratio $\frac{1}{\zeta(d-t)^2}  \prod^{d-1}_{j=d-t+1}\nolimits\frac{\beta_j}{R}$ is
\begin{equation*}
\left\{\begin{array}{ll} \frac{1}{R^{d-1}\prod^{d-1}_{j=2}n_j} , & t=d, \\ 
  \frac{1}{R^{t-1}\zeta(d-t)^2\prod^{d-1}_{j=d-t+1}n_j}, &1\leq t< d.
\end{array}\right.
\end{equation*}
\end{theorem}

Reducing to the best rank-1 approximation problem  $\max_{\norm{\mathbf u_j}=1,1\leq j\leq d} \innerprod{\mathcal A}{\bigotimes^d_{j=1}\mathbf u_j}$ (one can remove the square in the objective function),  i.e., $R=1,t=d$, from Theorem \ref{th:approx_bound} we have the following corollary. 
\begin{corollary}
	Let $\mathbf u_j$, $1\leq j\leq d$ be generated by Algorithm \ref{alg:main} with $R=1$ and $t=d$. Then the approximation ratio is $\frac{1}{\sqrt{\prod^{d-1}_{j=2}n_j  }}$ (either in deterministic or expected sense). 
	\end{corollary}

The above ratio is of the same order as \cite{he2010approximation,zhang2012cubic}. It is not the best to date \cite{so2010deterministic,he2014probability}; however, Algorithm \ref{alg:main} is more practical, and even in the context of rank-1 approximation problems,    Algorithm \ref{alg:main} with Procedures \ref{proc:get_v_from_M_B} and \ref{proc:get_v_from_M_C} is new   (we discuss it in the last paragraph of the next subsection).

\subsection{Discussions} \label{sec:approx_bound_discussions}
	\begin{figure}[h] 
	\centering
	\subfigure[Plots of the real ratio $\bigxiaokuohao{\sum^R_{i=1}\bignorm{\mathbf v_{j,i}}^2/\bigfnorm{M_{j,i}}^2}/R$  (in red) and the theoretical ratio $1/n$ (in black) for fourth-order tensors. $n$ varies from $10$ to $100$. Left: $j=1$; right: $j=2$.]{
		\label{fig:vm} %% label for second subfigure
		\includegraphics[height=3.5cm,width=13cm]{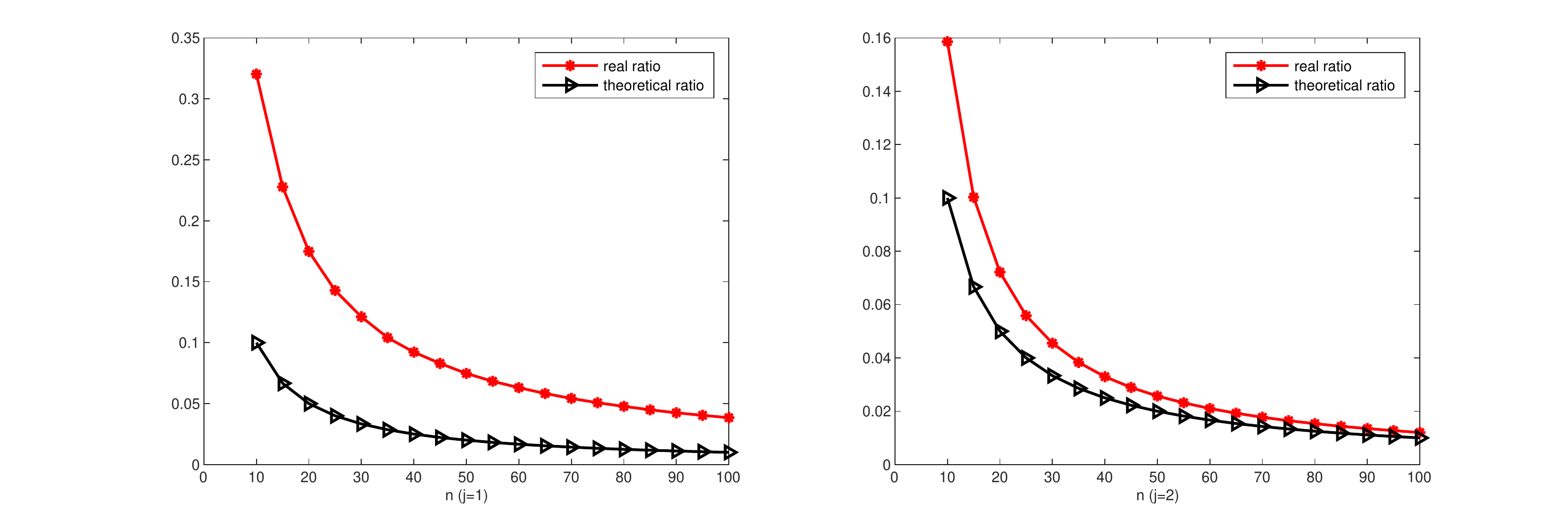}}\\
	\subfigure[Plots of the real ratio $\sum^R_{i=1}\innerprod{\mathbf u_{j,i}}{\mathbf v_{j,i}}^2/\bigfnorm{V_j}^2$ (in red) and the theoretical ratio $1/R$ (in black) where $\mathcal A\in\mathbb R^{100\times 100\times 100\times 100}$. $R$ varies from $10$ to $100$. Left: $j=1$; right: $j=2$.]{
		\label{fig:vuv} %% label for second subfigure
		\includegraphics[height=3.5cm,width=13cm]{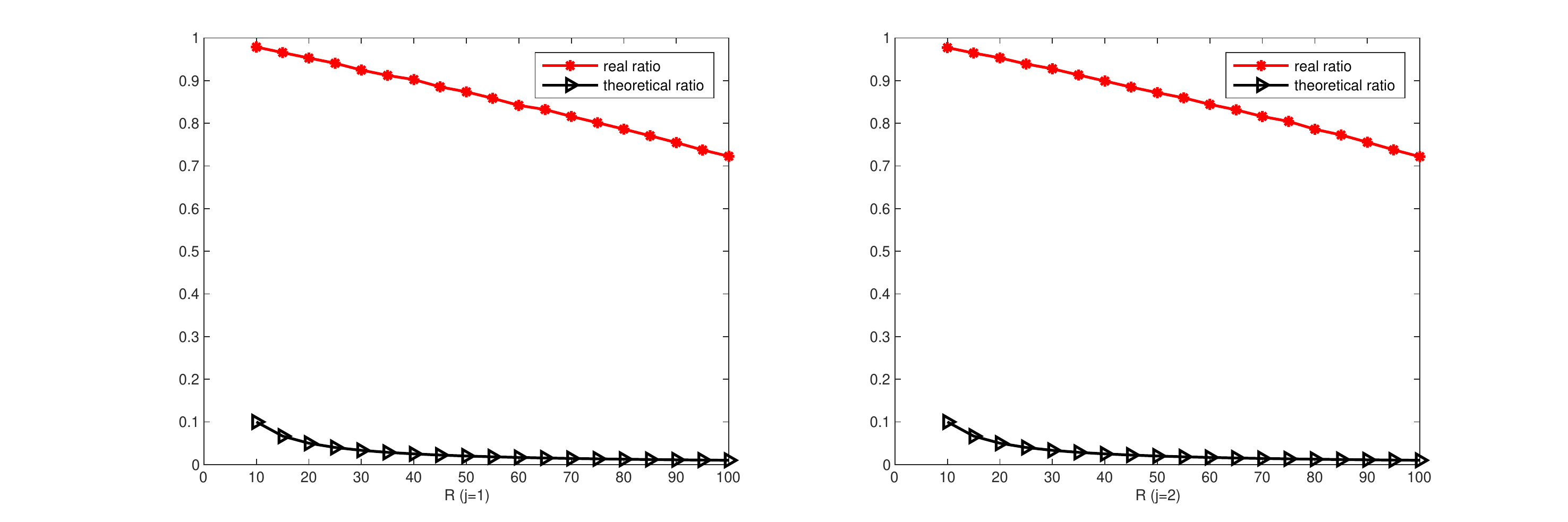}}
		\subfigure[Plots of real ratio versus theoretical ratio. Left: $\mathcal A\in\mathbb R^{n\times n\times n\times n}$ with $n$ varying from $10$ to $100$; right: $\mathcal A\in\mathbb R^{100\times 100\times 100\times 100} $ with $R$ varying from $10$ to $100$. In the figures: 
		Red: real ratio $\frac{G(U_1,\ldots,U_d)}{\sum^R_{i=1}\lambda_i(A_{(d)})^2}$; black with right arrow markers: $\frac{1}{n^2}$; black with diamond markers: $\frac{1}{R^3n^2}$.]{
		\label{fig:obj_up_lower_bound} %% label for second subfigure
		\includegraphics[height=3.5cm,width=13cm]{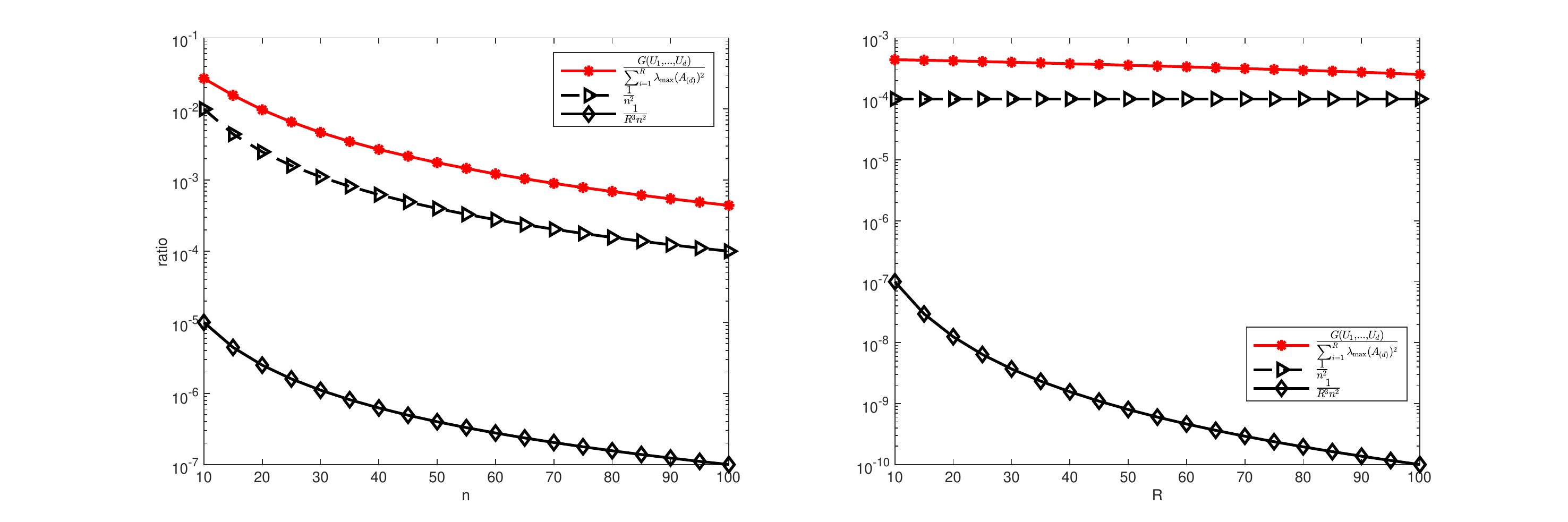}}
	\caption{Real ratio versus theoretical ratio.} 
	
	\label{fig:real_ratio_vs_th_ratio} %% label for entire figure
\end{figure} Several discussions concerning the approximation bound in Theorem \ref{th:approx_bound} are provided in this subsection.
 
	Firstly, Theorem \ref{th:approx_bound}   derives the worse-case theoretical approximation bound for general tensors. Looking into Lemma \ref{lem:approx_bound_informal_chain_inequality}, we see that the  approximation ratio  comes from two folds of each factor: the lower bound of $\bignorm{\mathbf v_{j,i}}^2/\bigfnorm{M_{j,i}}^2$ and the lower bound of $\sum^R_{i=1}\innerprod{\mathbf u_{j,i}}{\mathbf v_{j,i}}^2/\bigfnorm{V_j}^2$. The former reaches the worst bound $1/n_j$ when $M_{j,i}$  is row-wisely orthogonal, where each row has equal magnitude; the latter attains the worst bound $1/R$ when  $V_j$ is of rank-1 and each column admits equal magnitude.  These two worse cases rarely happen simultaneously, meaning that it is possible to further improve the theoretical lower bound. Currently, we cannot figure it out. Instead, we show the real ratio and the theoretical ratio via an example as follows.
	
	 We generate   tensors $\mathcal A\in\mathbb R^{n\times n\times n\times n}$, where every entry obeys the standard normal distribution. We set $t=4$ for the problem. $n$ varies from $10$ to $100$, where for each $n$, $50$ instances are generated. We use Procedure \ref{proc:get_v_from_M_A} in the algorithm. 
	In Fig. \ref{fig:vm}, we plot the curve  of the real ratio $\bigxiaokuohao{\sum^R\bignorm{\mathbf v_{j,i}}^2/\bigfnorm{M_{j,i}}^2}/R$ and curve of the theoretical ratio $1/n$ when $j=1,2$. 
		  The figure shows that although the real ratio is better,    it  still  approaches the theoretical one as $n$ increases. 
  In Fig. \ref{fig:vuv}, we fix the size of the tensor as $\mathcal A\in\mathbb R^{100\times 100\times 100\times 100}$, and vary $R$ from $10$ to $100$. We plot the real ratio $\sum^R_{i=1}\innerprod{\mathbf u_{j,i}}{\mathbf v_{j,i}}^2/\bigfnorm{V_j}^2$ and the theoretical one $1/R$. The figure   shows that the real ratio is far   better, while it also decreases as $R$ increases. 
  
  To further investigate the influence of $n$ and $R$ on the bound, we also plot the whole  real ratio $\frac{G(U_1,\ldots,U_d)}{\sum^R_{i=1}\lambda_{\max}(A_{(d)})^2}$ and the theoretical one $ {1}/{R^3n^2}$ in Fig. \ref{fig:obj_up_lower_bound}.  
  The tensors are the same as above, where in the left panel, we vary $n$ from $10$ to $100$ with $R=10$ fixed; in the right one, $R$ varies from $10$ to $100$ with $n=100$ fixed. The real ratio is in red, and the theoretical ratio is in black with diamond markers. We see that the real ratio decreased as $n$ increases, confirming the observation of Fig. \ref{fig:vm}; moreover, the real ratio is almost parallel to the theoretical ratio, showing that the theoretical ratio is reasonable up to a constant (if we   view $R$ as a constant). The right panel shows that when $R$ increases, the real ratio does not decrease as much as the theoretical one, showing that the term $1/R^3$ might be too loose in the theoretical ratio. Thus we also plot the curve of   $\frac{1}{n^2}$, which is in black with right arrow markers. We see that $1/n^2$ still under the curve of the real ratio, meaning that it is possible to improve the theoretical ratio in the absence of $R$. This still needs further research.

 Secondly, one may wonder whether   estimations such as $\bigfnorm{\mathcal A- \mathcal A^{\rm est}} \leq \alpha \bigfnorm{\mathcal A - \mathcal A^{\rm best}}$ can be established for the proposed algorithm for some factor $\alpha$, such as the truncated HOSVD \cite{de2000a}; here $\mathcal A^{\rm est}$ is constructed from $U^0_1,\ldots,U^0_d$, and $\mathcal A^{\rm best}$ denotes the optimal solution. This seems to be   difficult for a general tensor in the current CP format. The reason may be because unlike that every tensor can be written in the HOSVD format \cite[Theorem 2]{de2000a}, where the truncated HOSVD gives a good approximation, not every tensor admits a CP decomposition with orthonormal factors. Therefore, estimating $\bigfnorm{\mathcal A-\mathcal A^{\rm est}}$ is not easy for a general tensor.  
  
Thirdly, in the algorithm, other procedures for computing $\mathbf v_{j,i}$ from $M_{j,i}$ are possible. Another   randomized example is given in Appendix \ref{sec:auxiliary_procedures}. 
Lemma \ref{lem:approx_bound_get_v_from_M_D} shows that the factor $\beta_j$ of this procedure is worse than those of Procedures \ref{proc:get_v_from_M_A}-\ref{proc:get_v_from_M_C}. However, it is expected that this procedure might be more efficient.

 Finally, after this work is almost finished, we realize that in    problem \eqref{prob:ortho_main_max}, if     $R=1$, $t=d$,  and Procedure \ref{proc:get_v_from_M_A} is taken in Algorithm \ref{alg:main}, then the algorithm boils down essentially to  \cite[Algorithm 1 with DR2]{he2010approximation} for   best rank-1 approximation problems.   In this regard, Algorithm \ref{alg:main}  generalizes \cite[Algorithm 2]{he2010approximation} to orthogonal rank-$R$ approximation cases, and involves more procedures other than Procedure \ref{proc:get_v_from_M_A}. 

\subsection{Approximation results for nearly orthogonal tensors} We briefly consider the following type of tensors in this subsection:
\label{sec:approx_bound_nearly_orth}
\begin{assumption}\label{ass:structured_A}
\begin{equation}
\label{eq:structured_tensor}
\mathcal A = \sum^R_{i=1}\nolimits\sigma_i \bigotimesu + \mathcal E,
\end{equation}
where the last $t$ $U_j$'s are orthonormal, and the first $(d-t)$ $U_j$'s are columnwisely normalized. In addition, we assume that $\sigma_1>\sigma_2>\cdots >\sigma_R>0$. $\mathcal E$ denotes the noisy tensor. 
\end{assumption}

The cases that $\sigma_{i_1}=\sigma_{i_2}$ for some $1\leq i_1<i_2\leq R$ need more assumptions and analysis, which is out of the scope of this work, and will be studied in the future.

Running Algorithm \ref{alg:main} on $\mathcal A$ with $t\geq 2$, it is easy to obtain the following results, where for completeness, we provide a short discussion. Denote $U_{-d}\in\mathbb R^{\prod^{d-1}_{j=1}n_j\times R }$ as  the matrix whose  $i$-th column    is the vectorization of $\bigotimes^{d-1}_{j=1}\mathbf u_{j,i}$. 
\begin{proposition}
	\label{prop:structured_A_e0}
	Denote $(U^0_1,\ldots,U^0_d)$ as the factors generated by Algorithm \ref{alg:main}, where  the procedure \texttt{get\_v\_form\_M} takes  Procedures \ref{proc:get_v_from_M_A}, \ref{proc:get_v_from_M_B}, \ref{proc:get_v_from_M_C}.    If   Assumption \ref{alg:main} holds, where $t\geq 2$,  and we set the same $t$ in the algorithm;     in addition, if $\mathcal E = 0$, then
	\[
	G(U_1^0,\ldots,U^0_d) = \sum^R_{i=1}\nolimits\sigma_i^2 = G_{\max},~{\rm and}~U^0_j = U_j,~j=1,\ldots,d,
	\]
	where for $U^0_j = U_j$, we  ignore the sign, i.e., it in fact holds $\mathbf u_{j,i}^0 = \pm \mathbf u_{j,i}$ for each $j$ and $i$. 
\end{proposition}
\begin{proof}
	To achieve this, we first show that during the algorithm, 
	\begin{equation}\label{eq:proof:prop:structured_A_e0:1}
V_j = U_j\cdot \Sigma, ~j=d-t+1,\ldots,d,
	\end{equation}
where $\Sigma:= {\rm diag}(\sigma_1,\ldots,\sigma_R)$. 	If \eqref{eq:proof:prop:structured_A_e0:1} holds, then since $U_j^0$ is given by the polar decomposition of $V_j$, by Assumption \ref{ass:structured_A}, it follows
	 \begin{equation}\label{eq:proof:prop:structured_A_e0:2}
	 U^0_j=U_j, ~j=d-t+1,\ldots,d.
	 \end{equation}
	\eqref{eq:proof:prop:structured_A_e0:1} can be proved by induction method starting from $j=d$. When $j=d$, this is in fact given by   \eqref{eq:V_d}. To see this, we only need to show that
  $\lambda_i(A_{(d)}) = \sigma_i$ where we recall that $\lambda_i(\cdot)$ denotes the $i$-th largest singular value of a matrix.  Under Assumption \ref{ass:structured_A} when $t\geq 2$, $U_{-d}$ is orthonormal, giving that $  U_d\Sigma U^\top_{-d}$ is a reduced SVD of $A_{(d)}$, and the results follow.
    Assume that \eqref{eq:proof:prop:structured_A_e0:1} holds when $j=d,d-1,\ldots,m$; then $U^0_j=U_j,j=d,d-1,\ldots,m$.   Therefore, 
	$$\mathcal B_{m-1,i}=\mathcal A\times_m \mathbf u^{0\top}_{m,i}\times_{m+1}\times\cdots\times_d\mathbf u_{d,i}^{0\top} = \sigma_i \bigotimes^{m-1}_{j=1}\nolimits \mathbf u_{j,i},$$
	which is a rank-1 tensor. 
	Since $M_{j-1,i}$ is the unfolding of $\mathcal B_{m-1,i}$,  it can be seen that $\mathbf y_{m-1,i}$ generated by Procedures \ref{proc:get_v_from_M_A}, \ref{proc:get_v_from_M_B}, \ref{proc:get_v_from_M_C}  is just the vectorization of $\bigotimes^{m-2}_{j=1}\mathbf u_{j,i}$, and so    $\mathbf v_{m-1,i} = \sigma_i\mathbf u_{m-1,i}$, which demonstrates that \eqref{eq:proof:prop:structured_A_e0:1} holds when $j=m-1$. Thus \eqref{eq:proof:prop:structured_A_e0:1} and \eqref{eq:proof:prop:structured_A_e0:2} are valid for $j=d-t+1,\ldots,d$. 
	
	Since $\mathcal A\times_{d-t+1} \mathbf u^{0\top}_{d-t+1,i}\times_{m+1}\times\cdots\times_d\mathbf u_{d,i}^{0\top} $ is of rank-1, the rank-1 approximation procedure generates that $\mathbf u^0_{j,i}=\mathbf u_{j,i}$, $j=1,\ldots,d-t$. Hence, $G(U^0_1,\ldots,U^0_d) = G(U_1,\ldots,U_d) = \sum^R_{i=1}\sigma_i^2$. This completes the proof. 
\end{proof}

From the above discussions, and by    matrix perturbation theory, it is not hard to obtain the perturbed version of Proposition \ref{prop:structured_A_e0}. 
\begin{proposition}
	\label{prop:structured_A_enot0}
	Under the setting of Proposition \ref{prop:structured_A_e0}, where now $\mathcal E$ is   small enough,  we have $	G(U_1^0,\ldots,U^0_d) =   G_{\max} - O(\bignorm{\mathcal E})$,  where the constant behind the big O is nonnegative,  and
\[
\min\{\bignorm{\mathbf u_{j,i}^0+ \mathbf u_{j,i}},\bignorm{\mathbf u^0_{j,i}-\mathbf u_{j,i} }   \} = O(\bignorm{\mathcal E}),~j=1,\ldots,d,~i=1,\ldots,R.
\]
\end{proposition}

It remains to consider the case that $t=1$ in Assumption \ref{ass:structured_A}, i.e., only $U_d$ is orthonormal. In general, if there is no additional assumptions on $U_j$, $j=1,\ldots,d-1$, then it is hard to obtain the approximation results. For instance, if $A = PQ^\top$ where $P$ is orthonormal but $Q$ is not, then it is in general hard to recover $P$ and $Q$. We thus make the following incoherence assumption:
\begin{assumption}
	\label{ass:incoherence}
	There is at least a $U_j$, $1\leq j\leq d-1$ being incoherent with modules $0\leq \delta<1$, i.e.,
	\[
	\exists 1\leq j \leq d-1, ~{\rm with}~ \bigjueduizhi{\innerprod{\mathbf u_{j,i_1}}{\mathbf u_{j,i_2}}} \leq \delta,~\forall i_1\neq i_2.
	\]
\end{assumption}
It is clear that $U_j$ is a nearly orthonormal matrix if $\delta$ is small enough. In what follows, we assume without loss of generality that $U_{d-1}$ satisfies Assumption \ref{ass:incoherence}. %Let $U_{-d}\in\mathbb R^{\prod^{d-1}_{j=1}n_j \times R }$ still denote the matrix whose $i$-th column is the vectorization of $\bigotimes^{d-1}_{j=1}\mathbf u_{j,i}$.
We immediately have: 
\begin{proposition}
	If $U_{d-1}$ is incoherent with modules $\delta$, then $U_{-d}$ is also incoherent with modules $\delta$. 
\end{proposition}

We consider the expression of $A_{(d)} A^\top_{(d)} $:
\begin{small}
	\begin{eqnarray*}
A_{(d)}A^\top_{(d)} &=& U_d \Sigma U_{-d}^\top U_{-d} \Sigma U_d^\top + O(\mathcal E)\\
&=& U_d \Sigma^2 U_d^\top + U_d \Sigma   \left[\begin{smallmatrix}
	0& \mathbf u_{-d,1}^\top\mathbf u_{-d,2}  &\cdots & \mathbf u_{-d,1}^\top\mathbf u_{-d,R}\\ \vdots&0&\ddots & \vdots \\\vdots&\ddots&0 & \vdots \\ \mathbf u_{-d,R}^\top\mathbf u_{-d,1} & \cdots & \mathbf u_{-d,R}^\top\mathbf u_{-d,R-1} & 0
\end{smallmatrix} \right]    \Sigma U_d^\top + O(\mathcal E) \\
&=& U_d \Sigma^2 U_d^\top + O(\delta) + O(\mathcal E),
	\end{eqnarray*}
\end{small}
\!\!namely, $A_{(d)}A_{(d)}^\top$ is a perturbation of the eigen-decomposition $U_d \Sigma^2 U_d^\top $, given that $\delta$ and $\mathcal E$ are small enough.  If Assumption \ref{ass:structured_A} holds, by matrix perturbation theory, the above implies that 
\[
\min\{\bignorm{\mathbf u_{d,i}^0+ \mathbf u_{d,i}},\bignorm{\mathbf u^0_{d,i}-\mathbf u_{d,i} }   \} = O(\delta) +  O(\bignorm{\mathcal E}),~i=1,\ldots,R,
\]
where $U^0_d$ is generated by the algorithm. It then follows that       $ \mathcal A\times_d \mathbf u_{d,i}^{0\top}$ is a nearly rank-1 tensor with perturbation $O(\delta) + O(\mathcal E)$, and so 
$$\min\{\bignorm{\mathbf u_{j,i}^0+ \mathbf u_{j,i}},\bignorm{\mathbf u^0_{j,i}-\mathbf u_{j,i} }   \} = O(\delta) +  O(\bignorm{\mathcal E}),~j=1,\ldots,d-1.$$
We thus have a similar result as Proposition \ref{prop:structured_A_enot0}:
\begin{proposition}
	\label{prop:structure_A_enot0_t1}
		Denote $(U^0_1,\ldots,U^0_d)$ as the factors generated by Algorithm \ref{alg:main},  where   Assumption \ref{alg:main} holds with $t=1$; furthermore, Assumption \ref{ass:incoherence} holds.  If $\mathcal E$ and $\delta$ are small enough, and we also set $t=1$ in the algorithm, then $	G(U_1^0,\ldots,U^0_d) =   G_{\max} - O(\delta) - O(\bignorm{\mathcal E})$,  where the constants behind the big O are nonnegative,  and
		\[
\min\{\bignorm{\mathbf u_{j,i}^0+ \mathbf u_{j,i}},\bignorm{\mathbf u^0_{j,i}-\mathbf u_{j,i} }   \} = O(\delta) + O(\bignorm{\mathcal E}),~j=1,\ldots,d,~i=1,\ldots,R.		
		\]
\end{proposition}

	\section{Numerical Study}\label{sec:numer}
We evaluate the performance of Algorithm \ref{alg:main} with Procedures \ref{proc:get_v_from_M_A}, \ref{proc:get_v_from_M_B}, and \ref{proc:get_v_from_M_C}, and compare them with \cite[Procedure 3.1]{yang2019epsilon}.
All the   computations are conducted on an Intel i7-7770 CPU desktop computer with 32 GB of RAM. The supporting software is Matlab R2019b.  The Matlab  package Tensorlab  \cite{tensorlab2013} is employed for tensor operations.

\paragraph{Performance of   approximation algorithms}

      \begin{figure}[h] 
	\centering
	\subfigure[$t=2$. Left: $G(U_1,\ldots,U_d)$ versus $n$; right:   CPU time versus $n$.]{
		\label{fig:d3_iterative_l1} %% label for second subfigure
		\includegraphics[height=3.5cm,width=13cm]{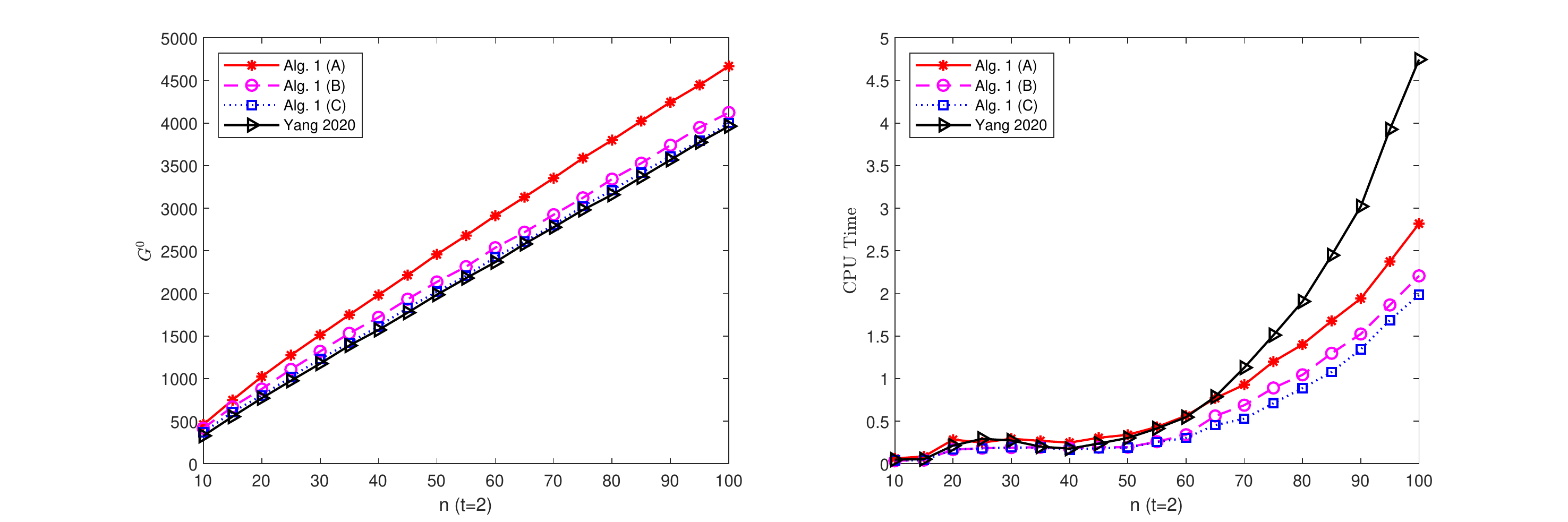}}\\
	\subfigure[$t=3$. Left: $G(U_1,\ldots,U_d)$ versus $n$; right:   CPU time versus $n$.]{
		\label{fig:d4_iterative_l1}%% label for second subfigure
		\includegraphics[height=3.5cm,width=13cm]{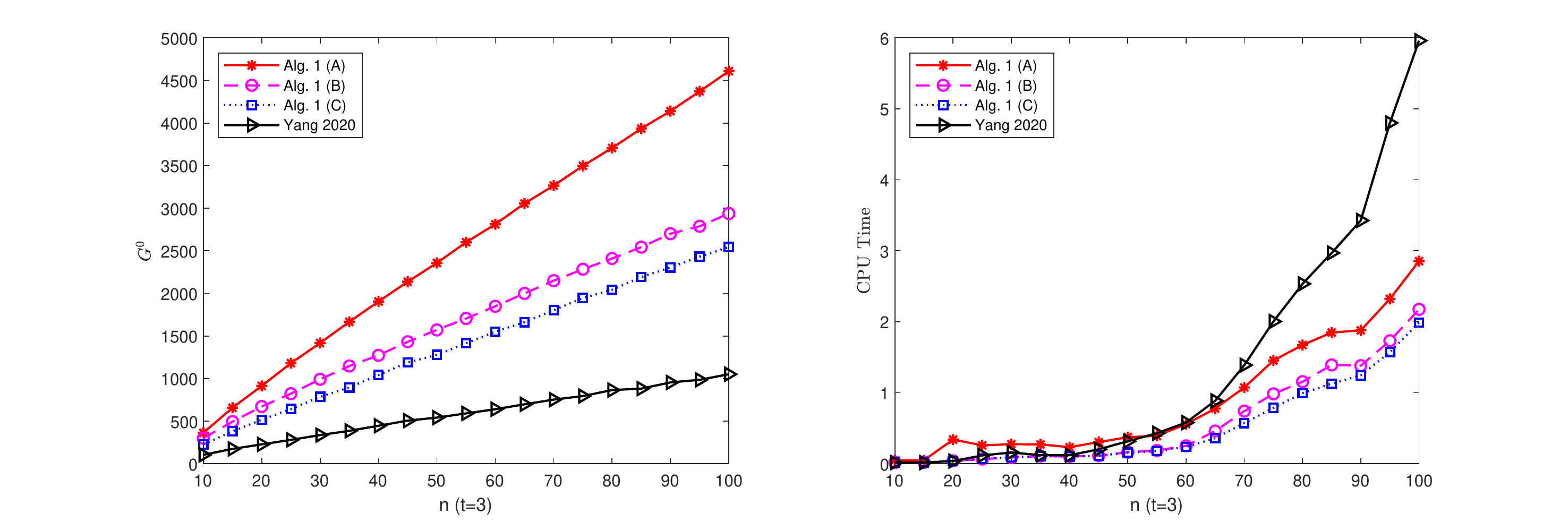}}\\
	\subfigure[$t=4$. Left: $G(U_1,\ldots,U_d)$ versus $n$; right:   CPU time versus $n$.]{
		\label{fig:d3d4_iterative_iterations_l1}%% label for second subfigure
		\includegraphics[height=3.5cm,width=13cm]{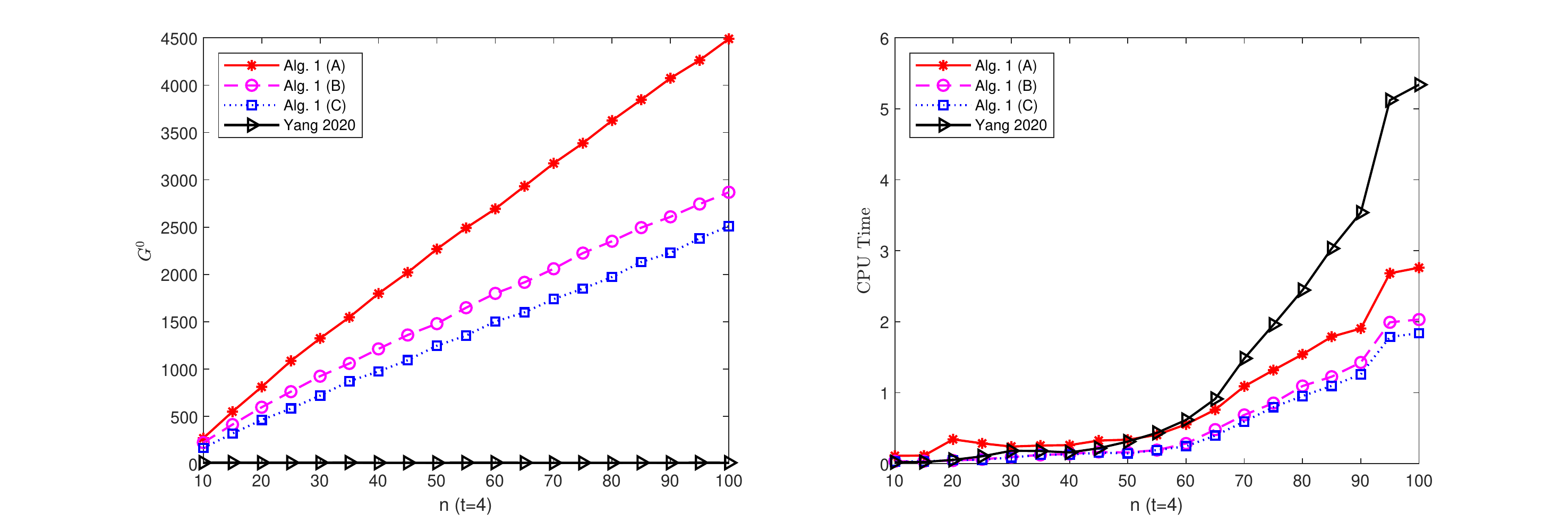}}\\
	\caption{Performances of  Alg. \ref{alg:main} \eqref{proc:get_v_from_M_A}  (red and star markers), Alg. \ref{alg:main} \eqref{proc:get_v_from_M_B} (magenta and circle markers), Alg. \ref{alg:main} \eqref{proc:get_v_from_M_C}  (blue and square markers), and \cite[Procedure 3.1]{yang2019epsilon} (black and rightarrow markers)  for approximately solving \eqref{prob:ortho_main_max} where $\mathcal A \in\mathbb R^{n\times n\times n\times n}$ are randomly generated. $n$ varies from $10$ to $100$.} 	
	\label{fig:randn_vary_n} %% label for entire figure
\end{figure}

We   compare Algorithm \ref{alg:main} with Procedures \ref{proc:get_v_from_M_A}, \ref{proc:get_v_from_M_B}, \ref{proc:get_v_from_M_C} with \cite[Procedure 3.1]{yang2019epsilon} in this part. They are respectively marked as Alg. \ref{alg:main} (\ref{proc:get_v_from_M_A}), Alg. \ref{alg:main} (\ref{proc:get_v_from_M_B}), Alg. \ref{alg:main} (\ref{proc:get_v_from_M_C}), and Yang2020. All the algorithms employ the same $\texttt{rank1approx}$ procedure  \cite[Procedure 3.2]{yang2019epsilon} for solving the rank-1 approximation subproblems. The tensors $\mathcal A \in\mathbb R^{n\times n\times n\times n}$ are randomly generated where every entry obeys the standard normal distribution. The presented results are averaged over $50$ instances for each case. 

We first fix $R=10$, and vary  $n$   from $10$ to $100$. The results are depicted in Fig. \ref{fig:randn_vary_n}. Alg. \ref{alg:main} (\ref{proc:get_v_from_M_A}) is in red with star markers, Alg. \ref{alg:main} (\ref{proc:get_v_from_M_B}) is in magenta with circle markers, Alg. \ref{alg:main} (\ref{proc:get_v_from_M_C}) is in blue with square markers, and Yang2020 is in black with rightarrow markers. The left panels show the curves of the objective values $G(U_1,\ldots,U_d)$ versus $n$, from which we observe that Alg. \ref{alg:main} (\ref{proc:get_v_from_M_A}) performs the best, followed by Alg. \ref{alg:main} (\ref{proc:get_v_from_M_B}); this is because that using SVD, Procedure \ref{proc:get_v_from_M_A} retains more information in $\mathbf v$ than other procedures. Yang2020 performs not as good as Algorithm \ref{alg:main}, and its performance gets worse when $t$ increases; this also has been pointed out in \cite[Table 3.1]{yang2019epsilon}, while the performance of Algorithm \ref{alg:main} is more stable. Therefore, we   conclude that Alg. \ref{alg:main} exploits the structure of the problem more than Yang2020.  
The right panels show the CPU time versus $n$, from which we see that Yang2020 needs more time than Alg. \ref{alg:main}, as Yang2020 requires to perform $t$ SVDs of size $n\times n^{d-1}$ (assuming $n_1=\cdots =n_d=n$), while the most expensive step of Alg. \ref{alg:main} is only   one SVD of this size. Considering Prodecures \ref{proc:get_v_from_M_A}, \ref{proc:get_v_from_M_B}, and \ref{proc:get_v_from_M_C}, there is no doubt that Alg. \ref{alg:main} (\ref{proc:get_v_from_M_A}) is the most expensive one among the three,  followed by Alg. \ref{alg:main} (\ref{proc:get_v_from_M_B}). The randomized version, i.e., Alg. \ref{alg:main} \eqref{proc:get_v_from_M_C}, is clearly the cheapest one.

     \begin{figure}[t] 
	\centering
		\includegraphics[height=3.5cm,width=13cm]{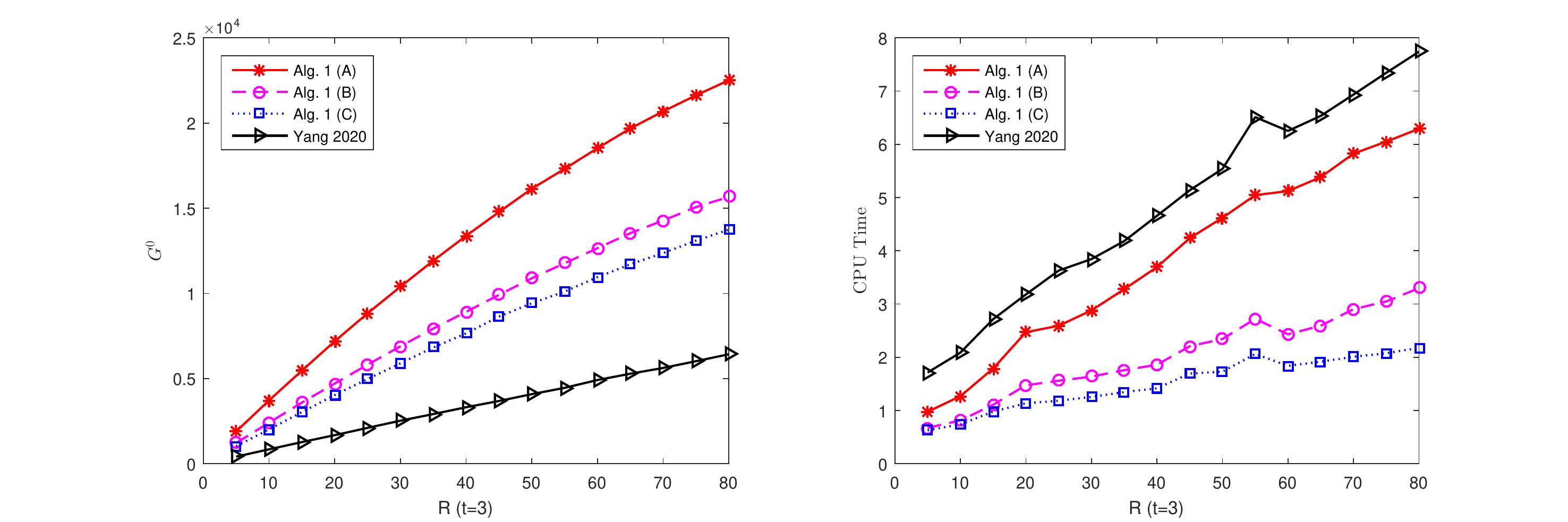}
  
	\caption{Performances of  Alg. \ref{alg:main} \eqref{proc:get_v_from_M_A}  (red and star markers), Alg. \ref{alg:main} \eqref{proc:get_v_from_M_B} (magenta and circle markers), Alg. \ref{alg:main} \eqref{proc:get_v_from_M_C}  (blue and square markers), and \cite[Procedure 3.1]{yang2019epsilon} (black and rightarrow markers)  for approximately solving \eqref{prob:ortho_main_max} where $\mathcal A \in\mathbb R^{80\times 80\times 80\times 80}$ are randomly generated. $R$ varies from $5$ to $80$.} 	
	\label{fig:randn_vary_r} %% label for entire figure
\end{figure}

We then fix $n=80$,  $t=3$,  and vary $R$ from $5$ to $80$. The results are plotted in Fig. \ref{fig:randn_vary_r}, from which we still can observe that Alg. \ref{alg:main} (\ref{proc:get_v_from_M_A}) performs the best in terms of the objective value. Considering the CPU time, Alg. \ref{alg:main} (\ref{proc:get_v_from_M_B}) and \ref{alg:main} (\ref{proc:get_v_from_M_C}) seem to be far more better; this shows the superiority of SVD-free procedures for computing the vector $\mathbf v$. 

\begin{table}[htbp]
	\renewcommand{\arraystretch}{1.9}
	\setlength{\tabcolsep}{6.5pt}
	\centering
	\caption{$\epsilon$-ALS initialized  by different strategies. $t=3$ and $R=10$ cases. `rel.err$^0$' and  `time$^0$' respectively represent  the relative error and CPU time evaluated at the initializers.}
	% Table generated by Excel2LaTeX from sheet 't=3'
	\begin{mytabular}{lllllll}
		\toprule
		\multicolumn{1}{l}{$n$} &       & \multicolumn{1}{c}{Yang2020}         & \multicolumn{1}{c}{Alg. \ref{alg:main} (\ref{proc:get_v_from_M_A})}        & \multicolumn{1}{c}{Alg. \ref{alg:main} (\ref{proc:get_v_from_M_B})}         & \multicolumn{1}{c}{Alg. \ref{alg:main} (\ref{proc:get_v_from_M_C})}         & \multicolumn{1}{c}{Random} \\
		\toprule
		\multirow{3}{*}{60}	& rel.err.(rel.err$^0$)	&   0.0139  	(0.0366) &	0.0140 	(0.0295) &	\textbf{0.0137} 	(0.0348) &	0.0218 	(0.0680) &	0.0291 \\
		& iter. & 7.32        & 6.26        & 8.70          & 8.04        &           18.72 \\
		& time(time$^0$)   & 0.59 (0.17)     & \textbf{0.52} (0.15)     &   0.56  (0.09)     &   0.53  (0.08)     &   0.92    \\
		\midrule
		\multirow{3}{*}{70}      & rel.err.(rel.err$^0$)	 & \textbf{0.0125}  (0.0255)  & \textbf{0.0125}    (0.0177)  &  0.0190    (0.0252)  & 0.0197    (0.0489)  & 0.0285             \\
		& iter. & 4.5      & 4.3         & 4.4        & 4.7              & 32.9  \\
		& time(time$^0$)   & 0.72 (0.36)    &  0.62  (0.26)    & \textbf{0.52} (0.16)    & 0.54 (0.16)    &   2.71   \\
		\midrule
		\multirow{3}{*}{80}      & rel.err.(rel.err$^0$)	 & 0.0191 (0.0239)  &  0.0186  (0.0216)  & 0.0189 (0.0218)  & \textbf{0.0114} (0.0343)  & 0.0193    \\
		& iter. & 19.2         & 9.2         & 9.7         & 24.4           & 35.08 \\
		& time(time$^0$)   & 3.16 (0.62)    &  1.77  (0.36)    & \textbf{1.73} (0.26)    & 3.50 (0.25)          & 4.58    \\
		\midrule
		\multirow{3}{*}{90}      & rel.err.(rel.err$^0$)	 & 0.0030 (0.0096)  & 0.0030 (0.0057)  & 0.0030 (0.0109)  & 0.0030 (0.0314)  & 0.0108   \\
		& iter. & 5.5       & 4.9        & 7.6        & 10.4  & 32.1  \\
		& time(time$^0$)   & 2.48 (1.02)    & \textbf{1.88} (0.54)    &  2.24  (0.40)    & 2.75 (0.36)    &   6.41   \\
		\bottomrule
	\end{mytabular}%
	\label{tab:t3}%
\end{table}%

\paragraph{Performance of approximation algorithms plus ALS for factor recovery}
      The tensors are generated similarly as \cite{yang2019epsilon,sorensen2012canonical}: $\mathcal A = \mathcal B/\bignorm{\mathcal B} + \beta\cdot\mathcal N/\|\mathcal N\|$, where $\mathcal B = \sum^R_{i=1} \sigma_i\bigotimesu $, $\mathcal N$ is an unstructured tensor, and $\beta$ denotes the noise level. Here we set $\beta=0.1$, $R=10$, $d=4$, and $n$ varies from $60$ to $90$.   $\sigma_i$,   $U_j$, and $\mathcal N$ are randomly drawn from a uniform distribution in $[-1,1]$. The last $t$ $U_j$ are  then orthonormalized by using the \texttt{orth} function, while the first $(d-t)$ ones are columnwisely normalized. We only test $t=3$ and $t=4$ cases. We compare $\epsilon$-ALS \cite{yang2019epsilon} with $\epsilon=10^{-8}$ initialized by different initializers generated respectively by Yang2020, Alg. \ref{alg:main} \eqref{proc:get_v_from_M_A}, Alg. \ref{alg:main} \eqref{proc:get_v_from_M_B}, and Alg. \ref{alg:main} \eqref{proc:get_v_from_M_C}. $\epsilon$-ALS initialized by random initializers is used as a baseline. 
The stopping criterion for $\epsilon$-ALS is 
$\sum^d_{j=1}\bigfnorm{U^{k+1}_j -U^k_j}/\bigfnorm{U^k_j}  \leq 10^{-5}$ 
or $k\geq 2000$.   The relative error is defined as follows \cite{sorensen2012canonical,yang2019epsilon}:
\vspace*{-0.2\baselineskip} 
\begin{equation*}\label{eq:relerr} \setlength\abovedisplayskip{3pt}
\setlength\abovedisplayshortskip{3pt}
\setlength\belowdisplayskip{3pt}
\setlength\belowdisplayshortskip{3pt}
{\rm rel.err} = \sum^d_{j=1}\nolimits\bignorm{U_{j} - U^{\rm out}_{j}\cdot\Pi_{j}}_F/\bignorm{U_{j} }_F,
\end{equation*}
where   $\Pi_{j}=\arg\min_{\Pi\in\boldsymbol{\Pi}} \bignorm{U_{j} - U^{\rm out}_{j}\cdot \Pi}_F$ and $\boldsymbol{\Pi}$ denotes the set of permutation matrices, and $U^{\rm out}_j$'s are the factors output by the iterative algorithm.  Finding the permutation can be efficiently done by the Hungarian algorithm \cite{kuhn1955hungarian}. The presented results are averaged over $50$ instances for each case.

\begin{table}[htbp]
	 	\renewcommand{\arraystretch}{1.5}
	\setlength{\tabcolsep}{6.5pt}
  \centering
	\caption{$\epsilon$-ALS initialized  by different strategies. $t=4$ and $R=10$ cases. `rel.err$^0$' and  `time$^0$' respectively represent  the relative error and CPU time evaluated at the initializers.}
    \begin{mytabular}{lllllll}
    	\toprule
\multicolumn{1}{l}{$n$} &       & \multicolumn{1}{c}{Yang2020}         & \multicolumn{1}{c}{Alg. \ref{alg:main} (\ref{proc:get_v_from_M_A})}        & \multicolumn{1}{c}{Alg. \ref{alg:main} (\ref{proc:get_v_from_M_B})}         & \multicolumn{1}{c}{Alg. \ref{alg:main} (\ref{proc:get_v_from_M_C})}         & \multicolumn{1}{c}{Random} \\
    \toprule
\multirow{3}{*}{60}      & rel.err.(rel.err$^0$)	 & 0.0211    (0.0328) & 0.0214    (0.0263)  & \textbf{0.0133}   (0.0320)  & 0.0209    (0.0658)  & 0.0706  \\
     & iter. & 13.4         & 12.0         & 11.3         & 18.6          & 36.7 \\
& time(time$^0$)   & 0.88 (0.19)   & 0.79 (0.16)    & {\bf0.69} (0.09)    & 1.01 (0.08)    & 1.75    \\
\midrule
\multirow{3}{*}{70} & rel.err.(rel.err$^0$)	 & 0.0129    (0.0300)  & 0.0129    (0.0209)  & \textbf{0.0126}   (0.0318)  & 0.0132    (0.0560)  & 0.0369  \\
      & iter. & 31.6          & 23.3          & 18.6          & 28.9         & 45.9 \\
& time(time$^0$)   & 3.16 (0.36)    & 2.21 (0.23)    & \textbf{1.76} (0.16)    &  2.57  (0.14)    & 4.03    \\
\midrule
\multirow{3}{*}{80} & rel.err.(rel.err$^0$)	 & 0.0279    (0.0347)  & \textbf{0.0276}   (0.0294)  & 0.0279   (0.0312)  & 0.0280    (0.0477)  & 0.0361  \\
      & iter. & 16.8          & 9.6          & 13.3         & 18.5          & 24.8 \\
& time(time$^0$)   & 2.96 (0.62)    & \textbf{1.80} (0.36)    &  2.19  (0.26)    & 2.78 (0.25)    & 3.30    \\
\midrule
\multirow{3}{*}{90} & rel.err.(rel.err$^0$)	 & 0.0029    (0.0064)  & 0.0029    (0.0041)  & 0.0029    (0.0079)  & 0.0029    (0.0414)  & 0.0194  \\
      & iter. & 30.6         & 11.5         & 20.1           & 16.9         & 43.5 \\
& time(time$^0$)   & 6.95 (0.98)   & \textbf{3.00} (0.50)    &  4.46  (0.38)    & 3.86 (0.35)    & 8.31    \\
\bottomrule
    \end{mytabular}%
  \label{tab:t4}%
\end{table}%

The results when $t=3$ and $t=4$ are reported in Tables \ref{tab:t3} and \ref{tab:t4}, in which `rel.err$^0$' and `time$^0$ respectively represent the relative error and CPU time evaluated at the initializers.  From the tables, we first observe that $\epsilon$-ALS initialized by Yang2020 and Alg. \ref{alg:main}  outperforms that with random initializations, considering both the relative error and computational time. Comparing  Alg. \ref{alg:main} with Yang2020, we can see that in most cases, $\epsilon$-ALS initialized by Alg. \ref{alg:main} performs comparable or slightly better. Comparing among   $\epsilon$-ALS initialized by the three versions of Alg. \ref{alg:main}, Alg. \ref{alg:main} \eqref{proc:get_v_from_M_C} is slightly worse in terms of the time; this is because although Alg. \ref{alg:main} \eqref{proc:get_v_from_M_C} is the most efficient, $\epsilon$-ALS initialized by it needs more iterations than the others. Therefore, it would be necessary to improve the efficiency of the iterative algorithms.

     \begin{figure}[t] 
	\centering
	\includegraphics[height=3.5cm,width=13cm]{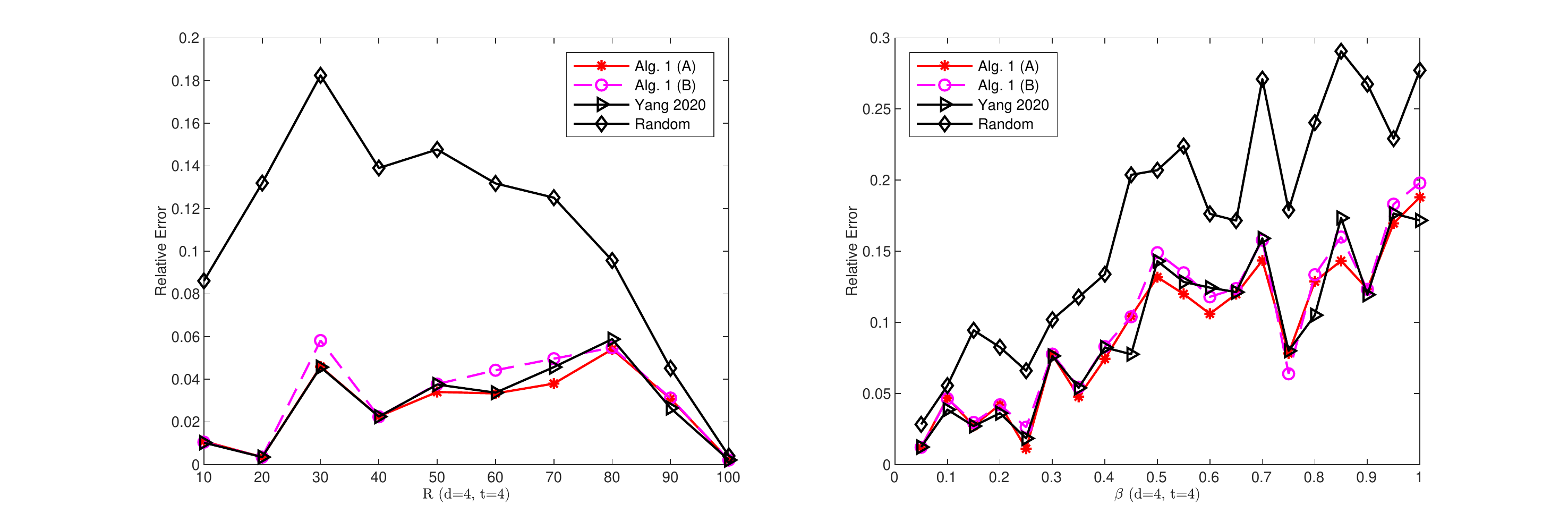}
	
	\caption{Comparisons of relative error  of $\epsilon$-ALS initialized by  Alg. \ref{alg:main} \eqref{proc:get_v_from_M_A}  (red and star markers), Alg. \ref{alg:main} \eqref{proc:get_v_from_M_B} (magenta and circle markers),  \cite[Procedure 3.1]{yang2019epsilon} (black and rightarrow markers), and random initializers (black and diamond markers). Left panel: $n$ is fixed to $100$, $\beta=0.1$, $t=4$, while $R$ varies from $10$ to $100$.  Right panel: $n$ is fixed to $50$, $t=4$, $R=10$, while the noise level $\beta$ varies from $0.05$ to $1$.} 	
	\label{fig:vary_R_vary_nl} %% label for entire figure
\end{figure}

Next, we fix $d=t=4$, $n=100$, $\beta=0.1$, and vary $R$ from $10$ to $100$; we also fix $n=50$, $R=10$, and vary $\beta$ from $0.05$ to $1$. Plots of relative error of different methods are depicted in Fig. \ref{fig:vary_R_vary_nl}, from which we still observe that $\epsilon$-ALS initialized by approximation algorithms performs better than that with random initializers, and the iterative algorithm initialized by Alg. \ref{alg:main} is comparable with that initialized by Yang2020.

\paragraph{CP approximation for clustering}
Tensor CP approximation for clustering works as follows: Suppose that we have $N$ samples $\mathcal A_1,\ldots,\mathcal A_n \in\T$, $d\geq 2$. For  given parameter $R$, and for unknown variables $A\in\mathbb R^{N\times R}, U_j\in\mathbb R^{n_j\times R}$, $j=1,\ldots,d$,  one solves the following problem first:
 \begin{equation}\label{prob:dim_reduction}  
\setlength\abovedisplayskip{4pt}
\setlength\abovedisplayshortskip{4pt}
\setlength\belowdisplayskip{4pt}
\setlength\belowdisplayshortskip{4pt}
\begin{split}
&\min  ~ \sum^N_{k=1}\bigfnorm{ \mathcal A_k - \sum^R_{i=1}\nolimits a^k_{i}\bigotimes^d_{j=1}\nolimits\mathbf u_{j,i}  }^2  \\
&~    {\rm s.t.}~~  A=(a^k_i) \in\mathbb R^{N\times R}, U_j^\top U_j = I,  1\leq j\leq d,
\end{split}
\end{equation}
where in $a^k_i$, $k$ represents the row index while $i$ is the column one. 
Write $\mathbf a^k:=[a^k_1,\ldots,a^k_R]^\top\in\mathbb R^R$, i.e., it is the transpose of the $k$-th row of $A$. The above problem means that one finds a common subspace in $\T$, which is spanned by the basis $\bigotimes^d_{j=1}\mathbf u_{j,i}$, $i=1,\ldots,R$, and projects the samples onto this subspace. The $\mathbf a^k$ can be seen as a representation of $\mathcal A_k$. This can be regarded as a dimension reduction procedure \cite{shashua2001linear,pesquet2001joint}. By stacking $\mathcal A_k$'s into a $(d+1)$-th order tensor $\mathcal A$ with $\mathcal A(k,:,\ldots,:)=\mathcal A_k$, and denoting $\mathbf a_i$ as the $i$-th column of $A$, the objective function of \eqref{prob:dim_reduction} can be written as $\bigfnorm{ \mathcal A -  \sum^R_{i=1}\mathbf a_i \otimes \mathbf u_{1,i}\otimes\cdots\otimes\mathbf u_{d,i} }^2$, and so \eqref{prob:dim_reduction} can be reduced to \eqref{prob:ortho_main_max}. 
Once $A=[\mathbf a_1,\ldots,\mathbf a_R]$ is obtained, we use $K$-means for clustering, with the transpose of   the rows of $A$, namely,  $\mathbf a^k$'s being the new samples. The cluster error is defined as follows: Let $K\geq 2$ be the given clusters,  $\phi_0:\T\rightarrow \{1,\ldots,K\}$ be the true clustering mapping, and $\phi$ be the estimated one. Denote ${\rm card}(\cdot)$ the cardinality of a set and $I(\cdot)$ the indicator function. Then \cite{sun2012regularized}
	\[
{\rm cluster~ err.} :=  \tbinom{N}{2}^{-1}{\rm card} \bigxiaokuohao{\{i,j\}: I\bigxiaokuohao{ \psi(\mathcal A_i)=\psi(\mathcal A_j)  }\neq I\bigxiaokuohao{\psi_0(\mathcal A_i)=\psi_0 (\mathcal A_j)},i<j} .
\]

We solve \eqref{prob:dim_reduction}  by $\epsilon$-ALS initialized by Alg. \ref{alg:main} \eqref{proc:get_v_from_M_A}, Alg. \ref{alg:main} \eqref{proc:get_v_from_M_B}, and by random initialization, with  $R=30$ for the problem.   For the first two initializations, the max iterations of $\epsilon$-ALS are set to $10$, while it is $1000$ for the random  initialization.    We also compare them with vanilla $K$-means, namely, the original samples $\mathcal A_k$'s are first vectorized and then clustered by $K$-means. The dataset
%Columbia Object Image Library 
COIL-100 \url{http://www.cs.columbia.edu/CAVE/software/softlib/coil-100.php} is used, which consists of $100$ objects, each containing $72$ images of size $128\times 128$ viewing from different angles. In our experiment, we each time randomly select $K=\{5,7,9,11,15,20 \}$ objects, each objects randomly selecting $M=50$ images, resulting into a third-order tensor $\mathcal A \in \mathbb R^{50K\times 128\times 128}$ that consists of $50K$ samples. For each case we run $50$ instances and present the averaged results in Table \ref{tab:clustering}.

 \begin{table}[htbp]
 	\renewcommand{\arraystretch}{2.8}
 	\setlength{\tabcolsep}{7pt}
 	\centering
 	\caption{CP approximation for clustering via first solving \eqref{prob:dim_reduction} by Alg. \ref{alg:main} \eqref{proc:get_v_from_M_A} + $\epsilon$-ALS, Alg. \ref{alg:main} \eqref{proc:get_v_from_M_B} + $\epsilon$-ALS, or random   + $\epsilon$-ALS; the iterative algorithm of the first two stops within $10$ iterations while that for the third is $1000$; then the $K$-means is performed to the reduced samples $\mathbf a^k$'s. We also use the vanilla $K$-means as the baseline. 
 	}
 	\begin{mytabular} {ccccccccc}       
 		\toprule
 		&  \multicolumn{2}{c}{Alg. \ref{alg:main} \eqref{proc:get_v_from_M_A}} &  \multicolumn{2}{c}{Alg. \ref{alg:main} \eqref{proc:get_v_from_M_B}} & \multicolumn{2}{c}{Random} &  \multicolumn{2}{c}{vanilla $K$-means}\\
 		\cmidrule(r){2-3} \cmidrule(r){4-5} \cmidrule(r){6-7} \cmidrule(r){8-9} 
 		$K$ & cluster err.  & time & cluster err.  & time  & cluster err. & time.& cluster err. & time \\
 		\toprule 
    5     & \textbf{1.10E-01} & 0.33  & 1.12E-01 & 0.24  & 1.44E-01 & 6.25  & 1.33E-01 & 0.44  \\
    7     & \textbf{9.60E-02} & 0.40  & 1.04E-01 & 0.31  & 1.64E-01 & 5.12  & 1.08E-01 & 0.79  \\
    9     & \textbf{8.24E-02} & 0.50  & 9.12E-02 & 0.42  & 1.40E-01 & 11.74  & 1.00E-01 & 1.25  \\
    11    & \textbf{7.48E-02} & 0.58  & 8.08E-02 & 0.50  & 1.21E-01 & 15.26  & 8.54E-02 & 1.76  \\
    15    & \textbf{6.31E-02} & 0.78  & 6.92E-02 & 0.68  & 9.56E-02 & 14.86  & 7.25E-02 & 3.23  \\
    20    & 5.36E-02 & 1.04  & \textbf{5.23E-02} & 0.98  & 7.84E-02 & 24.14  & 5.56E-02 & 5.58  \\
 		\bottomrule
 	\end{mytabular}
 	\label{tab:clustering}
 \end{table}

Considering the clustering error, we observe from the table that armed with the approximation algorithm, CP approximation based method achieves better performance than the vanilla $K$-means, while Alg. \ref{alg:main} \eqref{proc:get_v_from_M_A} is slightly better than Alg. \ref{alg:main} \eqref{proc:get_v_from_M_B}. However, if starting from a random initializer, CP approximation based method is worse than   vanilla $K$-means. This shows the usefulness of the introduced approximation algorithm. Considering the computational time, we see that the first two methods are also better than   vanilla $K$-means, which is because $\epsilon$-ALS is stopped within $10$ iterations, and because the sample size after reduction is $R=30$, while the sample size for the vanilla $K$-means is $128^2$. We also observe that Random + $\epsilon$-ALS usually cannot stop within $1000$ iterations, making it the slowest one. 
 
      \begin{figure}[h] 
 	\centering
 	\includegraphics[height=3.5cm,width=13cm]{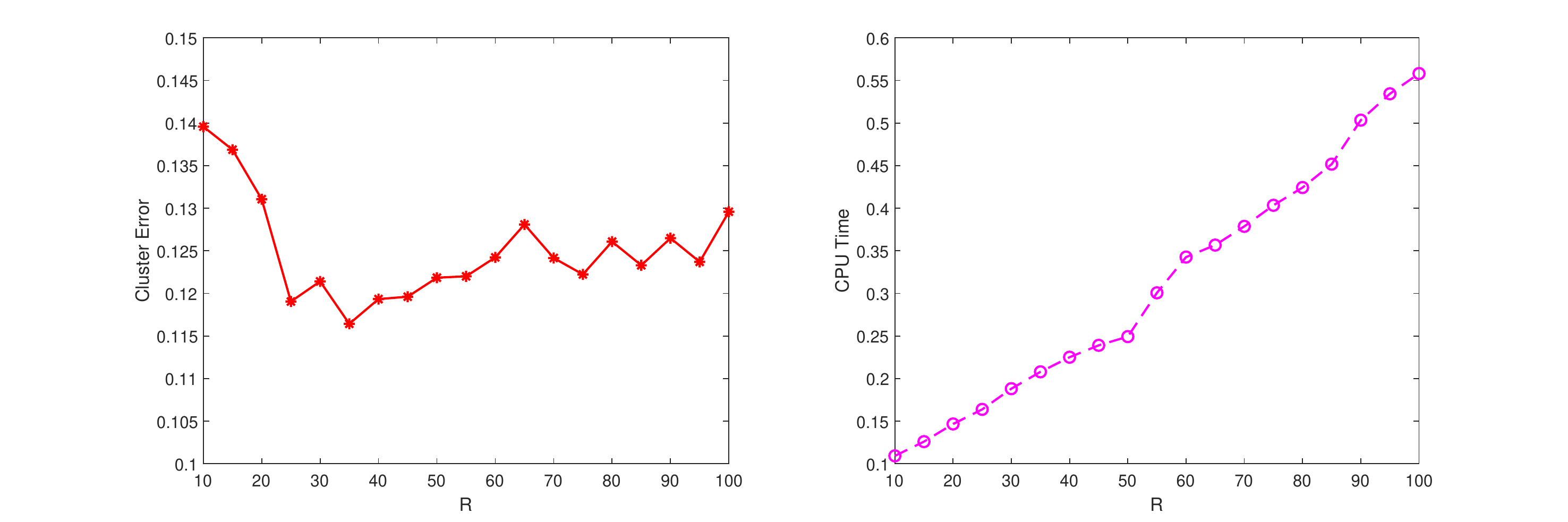}
 	
 	\caption{Cluster error and CPU time of  CP approximation for clustering of one instance with varying $R$ from $10$ to $100$, and $K=5$. Left: cluster error versus $R$; right: CPU time versus $R$.} 	
 	\label{fig:cluster_varying_R} %% label for entire figure
 \end{figure}

We next show the influence of $R$ on the performance. We fix an instance with $K=5$, and vary $R$ from $10$ to $100$. The cluster error and CPU time  of Alg. \ref{alg:main} \eqref{proc:get_v_from_M_A} + $\epsilon$-ALS is plotted in Fig. \ref{fig:cluster_varying_R}. We see that the cluster error does not change a lot when $R$ varies, while  around $R=30$, the cluster error seems to be slightly better than the other cases; this observation together with the CPU time shown in the right panel explains why we choose $R=30$ in our experiment. 

\section{Conclusions}\label{sec:conclusions}
 In \cite{yang2019epsilon}, an approximation procedure was proposed for  low-rank tensor approximation with orthonormal factors  by combining the truncated HOSVD and best rank-1 approximation. The  approximation bound  was only theoretically established when the number of latent orthonormal factors is one. To fill this gap, a modified approximation algorithm was developed in this work. It   allows either deterministic or randomized procedures to solve a key step of each latent orthonormal factor in the algorithm, giving some flexibilities.   The approximation bound, either in deterministic or expected sense, has been established  no matter how many orthonormal factors there are. 
 Moreover, compared with \cite[Procedure 3.1]{yang2019epsilon} which requires $t$ SVDs of size $n\times n^{d-1}$, the introduced approximation algorithm involves only one SVD of this size (plus other operations of smaller size), making it  more efficient. Numerical tests were provided to verify the usefulness of the algorithm, and its performance is favorably compared with    \cite[Procedure 3.1]{yang2019epsilon}.  The sequel work is to study approaches for finding   global solutions. A possible way is to use convex relaxation as those done for rank-1 approximation \cite{nie2014semidefinite,jiang2012tensor,yang2016rank}. However, extending the ideas from rank-1 approximation to rank-$R$ approximation is not so straightforward; this will be our next focus. 
 Another possible research thread is to extend the notion of best rank-1 approximation ratio of a tensor space \cite{qi2011best,li2018orthogonal}  to the best rank-$R$ approximation ratio setting (via problem \eqref{prob:ortho_main_max}) and study its properties.

 {\scriptsize\section*{Acknowledgement}  This work was supported by NSFC Grant 11801100  and the    Fok Ying Tong Education Foundation Grant 171094. We thank Mr. Xianpeng Mao for the help in the experiments.}

   \bibliography{../tensor,../TensorCompletion,../orth_tensor}
  \bibliographystyle{plain}

\appendix

\section{Another  $\texttt{get\_v\_from\_M}$ Procedure}\label{sec:auxiliary_procedures}

	%The final one is   randomized and very simple. A unit vector $\mathbf y$ is randomly and uniformly drawn from the unit sphere of $\mathbb R^m$, and then $\mathbf v = M\mathbf y$:
%
\begin{boxedminipage}{0.85\textwidth}\small
	\begin{equation}  \label{proc:get_v_from_M_D}
	\noindent {\rm Procedure}~ \mathbf v=  \texttt{get\_v\_from\_M}(M)  \tag{D}
	\end{equation}
	1. Randomly and uniformly draw a vector $\mathbf y\in\mathbb S^{m-1}$, where $\mathbb S^{m-1}$ denotes the unit sphere in $\mathbb R^m$;
	
	2. Return $\mathbf v = M\mathbf y$. 
\end{boxedminipage}
\begin{lemma}
	\label{lem:approx_bound_get_v_from_M_D}
	Let $\mathbf v\in\mathbb R^n$ be generated from $M\in\mathbb R^{n\times m}$ by Procedure \ref{proc:get_v_from_M_D}. Then it holds that $\E \bignorm{\mathbf v}^2 =   \frac{1}{m}\bigfnorm{V}^2$. 
	\end{lemma}   

\begin{proof}
	The proof is based on Lemmas \ref{lem:expectation_unit_sphere} and \ref{lem:expectation_inner_prod}. 
%		According to Lemma \ref{lem:expectation_inner_prod}, we have
%		\[
%		\E \bignorm{\mathbf v}^2 = \E \sum^{m}_{k=1} \nolimits \innerprod{M(k,:)}{\mathbf y}^2 = m^{-1}\sum^{m}_{k=1}\nolimits\bignorm{M(k,:)}^2 = m^{-1}\bigfnorm{M}^2.
%		\]
	\end{proof}
\begin{lemma}
	\label{lem:expectation_unit_sphere}
	Let $\mathbf y$ be randomly and uniformly drawn from the unit sphere $\mathbb S^{m-1}$ in $\mathbb R^m$. Then it holds that
	\[
	\E \mathbf (y^k)^2 =  {1}/{m},~\E  y^{k_1}  y^{k_2} = 0,
	\] 
	where $k,k_1,k_2=1,\ldots,m$, $k_1\neq k_2$, and $  y^k$ denotes the $k$-th entry of $\mathbf y$. 
	\end{lemma}
\begin{proof}
%	Since $\mathbf y$ is uniformly distributed on the unit sphere, by its symmetry, we have $\Einline \mathbf y(1)^2=\Einline\mathbf  y(2)^2 = \cdots = \Einline \mathbf y(n)^2$. This together with $\bignorm{\mathbf y}=1$ shows that $\Einline\mathbf y(k)^2 = \frac{1}{n}$ for each $k$. Such a property was also mentioned in \cite[p. 896]{he2014probability}. 
	The first property comes from  \cite[p. 896]{he2014probability}.  
	The second one   uses   symmetry. Let $\mathbf z\in\mathbb R^n$ be defined such that $  z^k=  y^k$ for all $k\neq k_1$, and $  z^{k_1}=-  y^{k_1}$. Then $\mathbf z$ is also uniformly distributed on the unit sphere, which means that $\Einline  z^{k_1}  z^{k_2} = \Einline   y^{k_1}  y^{k_2}$. By the definition of $\mathbf z$, this means that $\Einline  y^{k_1}  y^{k_2}=0$. Since $k_1$ and $k_2$ can be arbitrary, the results follow. 
	\end{proof}
\begin{lemma}
	\label{lem:expectation_inner_prod}
		Let $\mathbf y$ be randomly and uniformly drawn from   $\mathbb S^{m-1}$. For any constant vector $\mathbf a\in\mathbb R^m$,   it holds that
$
	\E \innerprod{\mathbf a}{\mathbf y}^2 =  m^{-1}{\bignorm{\mathbf a}^2}.
$
	\end{lemma}
\begin{proof}
	We have from Lemma \ref{lem:expectation_unit_sphere} that
 	\begin{eqnarray*}
	\E \innerprod{\mathbf a}{\mathbf y}^2 &=& \E \sum^m_{k_1=1}\nolimits\sum^m_{k_2=1}\nolimits  a^{k_1}  y^{k_1}  a^{k_2}  y^{k_2} \\
	&=& \sum^m_{k=1}\nolimits  (a^k)^2  (y^k)^2 + \sum_{k_1\neq k_2}\nolimits   a^{k_1}  a^{k_2}  y^{k_1}  y^{k_2} = m^{-1}{\bignorm{\mathbf a}^2} .
	\end{eqnarray*}
	\end{proof}

\end{document}